\pdfoutput=1

\makeatletter
\IfFileExists{./numapde-latex/numapde-packages.sty}{\providecommand*{\input@path}{}\edef\input@path{{./numapde-latex/}\input@path}}{}
\makeatother

\documentclass{numapde-preprint}

\usepackage{numapde-semantic}
\usepackage{numapde-style}
\usepackage{numapde-local}

\IfFileExists{./numapde-bibliography/numapde.bib}{\addbibresource{./numapde-bibliography/numapde.bib}}{\addbibresource{numapde.bib}}

\hypersetup{
	pdftitle={Dimensionally Consistent Preconditioning for Saddle-Point Problems},
	pdfauthor={Roland Herzog},
	pdfkeywords={saddle-point problems, preconditioning, dimensional consistency, physical units, minimum residual method}
}

\title{Dimensionally Consistent Preconditioning for Saddle-Point Problems}
\subtitle{}
\shorttitle{Dimensionally Consistent Preconditioning for Saddle-Point Problems}

\author{Roland Herzog\thanks{Technische Universität Chemnitz, Faculty of Mathematics, 09107 Chemnitz, Germany (\email{roland.herzog@mathematik.tu-chemnitz.de}, \url{https://www.tu-chemnitz.de/mathematik/part_dgl/people/herzog}, \orcid{0000-0003-2164-6575}).}}
\shortauthor{R. Herzog}

\dedication{}

\IfFileExists{numapde-uncolorizeHyperrefIfChanges.sty}{\RequirePackage{numapde-uncolorizeHyperrefIfChanges}}{}

\begin{document}
\maketitle

\begin{abstract}
The preconditioned iterative solution of large-scale saddle-point systems is of great importance in numerous application areas, many of them involving partial differential equations.
Robustness with respect to certain problem parameters is often a concern, and it can be addressed by identifying proper scalings of preconditioner building blocks.
In this paper, we consider a new perspective to finding effective and robust preconditioners.
Our approach is based on the consideration of the natural physical units underlying the respective saddle-point problem.
This point of view, which we refer to as dimensional consistency, suggests a natural combination of the parameters intrinsic to the problem.
It turns out that the scaling obtained in this way leads to robustness with respect to problem parameters in many relevant cases.
As a consequence, we advertise dimensional consistency based preconditioning as a new and systematic way to designing parameter robust preconditoners for saddle-point systems arising from models for physical phenomena.\end{abstract}



\section{Introduction}
\label{sec:introduction}

Saddle-point systems of the form 
\begin{equation}
	\label{eq:saddle-point_system}
	\begin{bmatrix}
		A & B_2^\star \\ B_1 & -C
	\end{bmatrix}
	\begin{pmatrix}
		u \\ p
	\end{pmatrix}
	= 
	\begin{pmatrix}
		f \\ g
	\end{pmatrix}
\end{equation}
arise in numerous applications, including computational fluid dynamics, the elastic deformation of solids, and quadratic programming.
We refer the reader to \cite{BenziGolubLiesen:2005:1} for a comprehensive treatment of these problems and their properties.
Often these systems arise from the discretization of partial differential equations and this will also be our focus and source of examples here.
Then \eqref{eq:saddle-point_system} is typically large-scale, and its efficient solution is of utmost importance in applications.

One method of choice are preconditioned iterative solvers of Krylov subspace type.
This has been an active area of research for several decades and we refer the reader to \cite{ElmanSilvesterWathen:2014:1,Wathen:2015:1,PestanaWathen:2015:1} and the references therein for recent surveys on various aspects of this topic.
The mesh independent convergence behavior is of particular importance.
One way to achieve it is to follow the principles of \emph{operator preconditioning}; see for instance \cite{Hiptmair:2006:1,MardalWinther:2011:1}.
Here one designs the preconditioner according to the mapping properties of the building blocks of the operator matrix on the left hand side of \eqref{eq:saddle-point_system}.

Another important aspect which has attracted a lot of attention recently is the issue of robust performance of preconditioned iterative solvers with respect to certain problem parameters; see for instance \cite{Klawonn:1998:2,SchoeberlZulehner:2007:1,MardalWinther:2011:1,Zulehner:2011:1,PearsonStollWathen:2012:1,KollmannZulehner:2013:1,ElvetunNielsen:2016:1,MardalNielsenNordaas:2016:1}. 
This requires that the blocks of the preconditioner scale appropriately with these parameters.
Exactly how robustness is achieved, however, seems to be a matter of experience and experiment.

In this paper, we propose another aspect which may help design and appropriately scale preconditioners for saddle-point problems but which seems to have gone largely unnoticed.
We propose to take into account the physical units of the primal and dual variables~$u$ and $p$, as well as the units of the first and second residuals.
For a clearer motivation, we shall restrict the discussion to self-adjoint systems with block-diagonal and self-adjoint positive definite preconditioners.
In this setting, the preconditioned minimum residual method (\minres) introduced in \cite{PaigeSaunders:1975:1} is the Krylov subspace method of choice; see also \cite{GuennelHerzogSachs:2014:1} and \cite[Chapter~4.1]{ElmanSilvesterWathen:2014:1}.
In order to monitor convergence, \minres evalutes the preconditioner-induced norm of the residual.
By respecting the physical units of the problem, our framework ensures that this residual norm is meaningful although both components of the residual in \eqref{eq:saddle-point_system} usually have completely different physical interpretations.

We mention that a proper choice of norm for minimum residual type methods to converge in has been the topic of numerous publications such as \cite{Wathen:2007:1,SilvesterSimoncini:2011:1,PestanaWathen:2015:1} and \cite[Chapter~4]{ElmanSilvesterWathen:2014:1}.
However we are not aware of a systematic treatment from the perspective of physical units.
We believe that this aspect supplements the idea of operator preconditioning, which takes into account the mapping properties of differential operators, and hope that it may facilitate the search for efficient and robust preconditioners in the future.

This paper is structured as follows.
\Cref{sec:general_framework} explains the framework of physically consistent preconditioning of self-adjoint saddle-point problems.
In \cref{sec:examples} we discuss a number of parameter dependent preconditioners, both from the literature and new ones, from the point of view of the underlying physical units.
Numerical results are presented in \cref{sec:numerical_experiments} which confirm the robustness of the preconditioners found, in each case with respect to \emph{all} problem parameters.
A conclusion concludes the paper.

A word on terminology is in order before we begin.
We distinguish physical \emph{dimensions} (such as length) from \emph{units} (such as meter, or \si{\meter}).
Strictly speaking, the technique introduced in this paper is one of dimensions.
Nevertheless, we prefer to work with physical units, which --- though mathematically equivalent --- we hope is more intuitive.
We will use the notation $\unit{u} = \si{\meter}$ to indicate that the unit associated with the variable~$u$ is meter.
We use SI units throughout but use the common abbreviations $\si{\newton} = \si{\kilo\gram\meter\per\second\squared}$ (Newton), $\si{\joule} = \si{\newton\meter}$ (Joule) and $\si{\watt} = \si{\joule\per\second}$ (Watt) where appropriate.

\section{A Dimensionally Consistent Framework for Saddle-Point Problems}
\label{sec:general_framework}

In this section we provide a framework for dimensionally consistent preconditioners, that is to say preconditioners respecting the underlying physical units of the respective problem.
As was mentioned in the introduction, we restrict the discussion to self-adjoint saddle-point problems
\begin{equation}
	\label{eq:self-adjoint_saddle-point_system}
	\begin{bmatrix}
		A & B^\star \\ B & -C
	\end{bmatrix}
	\begin{pmatrix}
		u \\ p
	\end{pmatrix}
	= 
	\begin{pmatrix}
		f \\ g
	\end{pmatrix}
	.
\end{equation}
The setting is made precise by means of the following example.

\subsection{Introductory Example: Stokes}
\label{subsect:introductory_example}
Let us begin with the discussion of the Stokes system, which describes slow viscous incompressible flows, here inside a bounded domain $\Omega \subset \R^3$.
The weak form of these equations, equipped for simplicity with no-slip boundary conditions, reads 
\begin{subequations}
	\label{eq:Stokes_system_with_bilinear_forms}
	\begin{alignat}{2}
		\label{eq:Stokes_system_with_bilinear_forms_1}
		a(u,v) + b(v,p) & = \dual{f}{v} & & \quad \text{for all } v \in V,
		\\
		\label{eq:Stokes_system_with_bilinear_forms_2}
		b(u,q) & = 0 & & \quad \text{for all } q \in Q,
	\end{alignat}
\end{subequations}
see for instance \cite[Chapter~3]{ElmanSilvesterWathen:2014:1}.
Here $u \in V = H^1_0(\Omega)^3$ denotes the primal variable (velocity) and $p \in Q = L^2(\Omega)/\R \simeq L^2_0(\Omega)$ is the pressure or dual variable.
Here $L^2_0(\Omega)$ denotes the space of $L^2(\Omega)$ functions with zero mean.
The bilinear forms are given by 
\begin{equation}
	\label{eq:Stokes_forms}
	a(u, v) = \mu \int_\Omega \nabla u \dprod \nabla v \dx, \quad
	b(u, q) = - \int_\Omega q \, \div u \dx, \quad
	\dual{f}{v} = \int_\Omega F \cdot v \dx.
\end{equation}
The constant $\mu$ denotes the dynamic viscosity of the fluid under consideration and $F$ is a given volume force density describing, \eg, the influence of gravity.
As usual, we can switch between the variational form \eqref{eq:Stokes_system_with_bilinear_forms} and the operator form \eqref{eq:self-adjoint_saddle-point_system} by letting $A \in \cL(V,V^*)$ and $B \in \cL(V,Q^*)$ be defined according to
\begin{equation*}
	\dual{Au}{v} \coloneqq a(u,v), \quad
	\dual{Bu}{q} \coloneqq b(u,q).
\end{equation*}
Here $V^*$ and $Q^*$ are the dual spaces of $V$ and $Q$, respectively, and $\dual{\cdot}{\cdot}$ denotes the duality pairing.
Moreover, the adjoint operator $B^\star \in \cL(Q,V^*)$ is given by $\dual{B^\star p}{v} = b(v,p)$ and $C = 0$ holds for the Stokes example.

Associated with problem \eqref{eq:Stokes_system_with_bilinear_forms} is the Lagrangian
\begin{equation}
	\label{eq:Stokes_Lagrangian}
	\begin{aligned}
		\cL(u,p) 
		&
		= 
		\frac{1}{2} a(u,u) + b(u,p) - \dual{f}{u}
		\\
		&
		=
		\frac{1}{2} \mu \int_\Omega \nabla u \dprod \nabla u \dx - \int_\Omega p \, \div u \dx - \int_\Omega F \cdot u \dx
		.
	\end{aligned}
\end{equation}
Let us consider its physical unit and recall that the unit of velocity is $\unit{u} = \si{\meter\per\second}$.
Moreover, we have $\unit{p} = \si{\newton\per\meter\squared}$, $\unit{F} = \si{\newton\per\meter\cubed}$ and $\unit{\mu} = \si{\newton\second\per\meter\squared}$.
Recall moreover that integration over $\Omega \subset \R^3$ adds a factor of $\si{\meter\cubed}$ while the differentiation with respect to spatial coordinates (by $\nabla$ and $\div$) adds a factor of $\si{\per\meter}$. 
Consequently we obtain the following units for each of the three terms in \eqref{eq:Stokes_Lagrangian}:
\sisetup{per-mode = fraction}
\begin{equation*}
	\begin{aligned}
		\unit[auto]{\mu \int_\Omega \nabla u \dprod \nabla u \dx} 
		&
		=
		\si{\newton\second\per\meter\squared} \cdot \si{\meter\cubed} \cdot \si{\per\meter} \cdot \si{\meter\per\second} \cdot \si{\per\meter} \cdot \si{\meter\per\second}
		=
		\si{\newton\meter\per\second}
		=
		\si{\watt}
		,
		\\
		\unit[auto]{\int_\Omega p \, \div u \dx} 
		&
		=
		\si{\meter\cubed} \cdot \si{\newton\per\meter\squared} \cdot \si{\per\meter} \cdot \si{\meter\per\second}
		=
		\si{\newton\meter\per\second}
		=
		\si{\watt}
		,
		\\
		\unit[auto]{\int_\Omega F \cdot v \dx} 
		&
		=
		\si{\meter\cubed} \cdot \si{\newton\per\meter\cubed} \cdot \si{\meter\per\second}
		=
		\si{\newton\meter\per\second}
		=
		\si{\watt}
		.
	\end{aligned}
\end{equation*}
The unit of the Lagrangian is thus $\si{\watt}$ (Watt).

The residual $r_1 \in V^*$ associated with \eqref{eq:Stokes_system_with_bilinear_forms_1}, is defined by $\dual{r_1}{v} \coloneqq - \cL_u(u,p) \, v = \dual{f}{v} - a(u,v) - b(v,p)$.
Clearly, $\unit{\dual{r_1}{v}} = \si{\watt}$ holds and since units are multiplicative, we obtain
\begin{equation*}
	\unit{r_1} 
	=
	\frac{\si{\watt}}{\unit{v}} 
	=
	\si{\watt\second\per\meter}
	=
	\si{\newton}
	.
\end{equation*}
We are thus reminded that the first equation \eqref{eq:Stokes_system_with_bilinear_forms_1}, is a balance of forces, measured in $\si{\newton}$ (Newton).
Similarly, we obtain for the second residual $r_2 \in Q^*$ defined by $\dual{r_2}{q} \coloneqq - \cL_p(u,p) \, q = - b(u,q)$ that
\begin{equation*}
	\unit{r_2} 
	=
	\frac{\si{\watt}}{\unit{q}} 
	=
	\si{\watt\meter\squared\per\newton}
	=
	\si{\meter\cubed\per\second}
	.
\end{equation*}
This reminds us that the second equation \eqref{eq:Stokes_system_with_bilinear_forms_2} is a balance of volumetric fluid production rates.

Summarizing the findings so far, we are led to endow the spaces $V$ and $Q$ for the primal and dual variables as well as their duals $V^*$ and $Q^*$ for the components of the residual with the following units:
\begin{equation*}
	\begin{aligned}
		\unit{V} & = \si{\meter\per\second},
						 &
		\unit{V^*} & = \si{\newton},
		\\
		\unit{Q} & = \si{\newton\per\meter\squared},
						 &
		\unit{Q^*} & = \si{\meter\cubed\per\second}.
	\end{aligned}
\end{equation*}
The product of the units in each row equals the unit of the Lagrangian, $\unit{\cL(u,p)} = \si{\newton\meter\per\second} = \si{\watt}$.

We proceed to discuss the role of the preconditioner.
As was mentioned in the introduction, we restrict the discussion to block-diagonal preconditioners,
\begin{equation*}
	P = 
	\begin{bmatrix}
		P_V & 0 \\ 0 & P_Q
	\end{bmatrix}
	.
\end{equation*}
Moreover, we assume that $P_V \in \cL(V,V^*)$ and $P_Q \in \cL(Q,Q^*)$ are self-adjoint and positive definite, as required by the preconditioned minimum residual method (\minres).
We observe that $P_V$ induces an inner product and hence a norm on $V$ by virtue of $(u,v)_{P_V} \coloneqq \dual{P_V u}{v}$ and $\norm{u}_{P_V}^2 = \dual{P_V u}{u}$.
Consequently, it also induces an inner product and norm on the dual space, namely $(r,s)_{P_V^{-1}} \coloneqq \dual{r}{P_V^{-1} s}$ and $\norm{r}_{P_V^{-1}}^2 = \dual{r}{P_V^{-1} r}$.
Similar considerations apply to $P_Q$ and the spaces $Q$ and $Q^*$.

Notice that \minres monitors the (squared) norm of the residual in the preconditioner-induced norm, i.e.,
\begin{equation}
	\label{eq:Stokes_squared_residual_norm}
	\norm{r}_{P^{-1}}^2
	=
	\norm{r_1}_{P_V^{-1}}^2
	+
	\norm{r_2}_{P_Q^{-1}}^2
	=
	\dual{r_1}{P_V^{-1} r_1}
	+ 
	\dual{r_2}{P_Q^{-1} r_2}
	.
\end{equation}
From this we observe that the role of the preconditioner is not only to provide a mapping from $V \times Q$ onto $V^* \times Q^*$ which respects the underlying Sobolev spaces (the perspective of operator preconditioning), but also to take into account the appropriate physical units assigned to these spaces.
Only when this is the case, is the quantity \minres computes in \eqref{eq:Stokes_squared_residual_norm} physically meaningful.
In our present Stokes example, we have $\unit[big]{\norm{r}_{P^{-1}}^2} = \si{\watt}$.

Let us confirm that a frequently used preconditioner for \eqref{eq:Stokes_system_with_bilinear_forms} does respect the physical units.
Our preconditioner of choice is one with diagonal blocks
\begin{equation}
	\label{eq:Stokes_preconditioner_blocks}
	\begin{aligned}
		\dual{P_V u}{v} 
		& 
		=
		\mu \int_\Omega \nabla u \dprod \nabla v \dx
		=
		a(u,v)
		,
		\\
		\dual{P_Q p}{q}
		&
		=
		\mu^{-1} \int_\Omega p \, q \dx
		,
	\end{aligned}
\end{equation}
see for instance \cite{WathenSilvester:1993:1} and \cite[Chapter~4]{ElmanSilvesterWathen:2014:1}.
In fact in these references the Stokes problem was considered in a dimensionless form, i.e., $\mu$ was replaced by one.
However the scaling \eqref{eq:Stokes_preconditioner_blocks} was proposed in \cite{UrRehmanGeenenVuikSegalMacLachlan:2011:1}.
Clearly, \eqref{eq:Stokes_preconditioner_blocks} is only an \emph{ideal} preconditioner which is too costly to realize in practice.
However the subsequent analysis is not affected when $P$ is replaced by a spectrally equivalent operator such as a geometric multigrid scheme.

In order to verify that $P$ properly handles the physical units, we only need to confirm that $\unit{P_V u} = \unit{V^*}$ and $\unit{P_Q p} = \unit{Q^*}$ holds.
Indeed, we obtain
\begin{equation*}
	\begin{aligned}
		\unit{P_V u}
		&
		=
		\unit[auto]{\mu \int_\Omega \nabla u \dprod \nabla \cdot \dx}
		=
		\si{\newton\second\per\meter\squared} \cdot \si{\meter\cubed} \cdot \si{\per\meter} \cdot \si{\meter\per\second} \cdot \si{\per\meter} 
		=
		\si{\newton}
		=
		\unit{V^*}
		,
		\\
		\unit{P_Q p}
		&
		=
		\unit[auto]{\mu^{-1} \int_\Omega p \, \cdot \dx}
		=
		\si{\meter\squared\per\newton\per\second} \cdot \si{\meter\cubed} \cdot \si{\newton\per\meter\squared} 
		=
		\si{\meter\cubed\per\second}
		=
		\unit{Q^*}
		.
	\end{aligned}
\end{equation*}
Moreover, it is easy to confirm that the relevant constants
\begin{equation}
	\label{eq:error_estimate_constants}
	\begin{aligned}
		\norm{A} 
		&
		\coloneqq 
		\sup_{u \in V \setminus \{0\}} \frac{\dual{A u}{u}}{\norm{u}_{P_V}^2},
		&
		\alpha 
		&
		\coloneqq 
		\inf_{u \in \ker B \setminus\{0\}} \frac{\dual{A u}{u}}{\norm{u}_{P_V}^2},
		\\
		\norm{B}
		&
		\coloneqq 
		\sup_{p \in Q \setminus \{0\}} \sup_{u \in V \setminus \{0\}} \frac{\dual{Bu}{p}}{\norm{u}_{P_V} \norm{p}_{P_Q}}
		&
		\beta 
		&
		\coloneqq 
		\inf_{p \in Q \setminus \{0\}} \sup_{u \in V \setminus \{0\}} \frac{\dual{Bu}{p}}{\norm{u}_{P_V} \norm{p}_{P_Q}}
	\end{aligned}
\end{equation}
are all independent of $\mu$ and, indeed, dimensionless.
Consequently, \minres verifies a convergence bound which is uniform in $\mu$.
We can thus conclude that the scaling of the preconditioner blocks as in \eqref{eq:Stokes_preconditioner_blocks} achieves both robustness and dimensional consistency, i.e., physical significance, at the same time.
This observation is indicative for all examples throughout the paper and it motivates our proposal to consider the physical units in search for parameter robust preconditioners.

\subsection{General Setting}
\label{subsec:general_setting}

Consider the saddle-point problem \eqref{eq:self-adjoint_saddle-point_system}.
Suppose that $V$ and $Q$ are Hilbert spaces, each of which in addition bears a physical unit.
We point out that these physical units of the primal and dual variables are known from the modelling.
We assume that $A \in \cL(V,V^*)$, $B \in \cL(V,Q^*)$ and $C \in \cL(Q,Q^*)$ as well as $f \in V^*$ and $g \in Q^*$ are bounded linear operators.
We further assume that the Lagrangian associated with \eqref{eq:self-adjoint_saddle-point_system}, 
\begin{equation}
	\label{eq:Lagrangian}
	\begin{aligned}
		\cL(u,p) 
		&
		= 
		\frac{1}{2} \dual{Au}{u} + \dual{Bu}{p} - \frac{1}{2} \dual{Cp}{p} - \dual{f}{u} - \dual{g}{p}
		,
	\end{aligned}
\end{equation}
is dimensionally consistent, i.e., that all terms in \eqref{eq:Lagrangian} bear the same physical unit $\unit{\cL(u,p)}$.
We then equip the dual spaces with physical units according to
\begin{equation*}
	\unit{V^*} \coloneqq \frac{\unit{\cL(u,p)}}{\unit{V}}, \quad
	\unit{Q^*} \coloneqq \frac{\unit{\cL(u,p)}}{\unit{Q}}.
\end{equation*}

Consider now a block-diagonal preconditioner 
\begin{equation}
	\label{eq:preconditioner}
	P = 
	\begin{bmatrix}
		P_V & 0 \\ 0 & P_Q
	\end{bmatrix}
	,
\end{equation}
where $P_V \in \cL(V,V^*)$ and $P_Q \in \cL(Q,Q^*)$ are self-adjoint and positive definite (spd).
We call the preconditioner $P$ \emph{dimensionally consistent} if $\unit{P_V u} = \unit{V^*}$ and $\unit{P_Q p} = \unit{Q^*}$ holds.
This is equivalent to the condition $\unit{\dual{P_V u}{u}} = \unit{\cL(u,p)} = \unit{\dual{P_Q p}{p}}$.
Notice that a dimensionally consistent preconditioner renders the unit of the squared residual norm \eqref{eq:Stokes_squared_residual_norm} well-defined, which then equals the unit of the Lagrangian.

\section{Examples}
\label{sec:examples}

In this section we discuss a number of saddle-point problems and corresponding block-diagonal, self-adjoint and positive definite preconditioners to be used with \minres.
In each example, the dimensional consistency of the preconditioner leads to its robustness w.r.t.\ \emph{all} parameters of the respective problem, in addition to mesh independence.

\subsection{Nearly Incompressible Elasticity}
\label{subsec:elasticity}

We consider the deformation of a body $\Omega \subset \R^3$ belonging to the class of linear, isotropic elastic materials with an emphasis on the nearly incompressible limit.
As is customary for these materials, we introduce an extra variable $p$ for the hydrostatic pressure (see for instance \cite{Wieners:2000:1}) in order to overcome the ill-conditioning of a purely displacement-based formulation known as \emph{locking} \cite{BabuskaSuri:1992:1}.
We employ the standard isotropic stress-strain relation, $\sigma = 2 \, \mu \, \varepsilon(u) + p \, I$ and $\div u = \lambda^{-1} p$.
Here $\varepsilon(u) = (\nabla u + \nabla u^\top)/2$ denotes the symmetrized Jacobian of $u$, while $\mu$ and $\lambda$ denote the Lamé constants.
The nearly incompressible case is obtained when $\lambda \gg \mu$.

We consider a problem where the deformed body is clamped on part of the domain boundary $\Gamma_D$ while traction forces $F$ act on the remaining part $\Gamma_N$.
The variational mixed formulation obtained in this way is described by the spaces $V = \{ v \in H^1(\Omega;\R^3): v = 0 \text{ on } \Gamma_D \}$ for the displacement and $Q = L^2(\Omega)$ for the hydrostatic pressure.
The bilinear and linear forms associated with this problem are
\begin{equation}
	\label{eq:elasticity_forms}
	\begin{aligned}
		a(u, v) 
		& 
		= 
		2 \, \mu \int_\Omega \varepsilon(u) \dprod \varepsilon (v) \dx, 
		&
		b(u, q) 
		& 
		= 
		\int_\Omega q \, \div u \dx, 
		\\
		c(p,q) 
		&
		= 
		\lambda^{-1} \int_\Omega p \, q \dx,
		&
		f(v) 
		&
		= 
		\int_{\Gamma_N} F \cdot v \dx
	\end{aligned}
\end{equation}
as well as $g = 0$.
The physical units associated with the spaces $V$ and $Q$ and their duals are
\begin{equation*}
	\begin{aligned}
		\unit{V} & = \si{\meter},
						 &
		\unit{V^*} & = \si{\newton},
		\\
		\unit{Q} & = \si{\newton\per\meter\squared},
						 &
		\unit{Q^*} & = \si{\meter\cubed},
	\end{aligned}
\end{equation*}
and the units of the remaining data is $\unit{\mu} = \unit{\lambda} = \unit{F} = \si{\newton\per\meter\squared}$.
The unit of the Lagrangian \eqref{eq:Lagrangian} is thus $\unit{\cL(u,p)} = \si{\newton\meter} = \si{\joule}$ (Joule).

We consider the preconditioner \eqref{eq:preconditioner} with blocks
\begin{equation}
	\label{eq:elasticity_preconditioner_blocks}
	\begin{aligned}
		\dual{P_V u}{v} 
		& 
		=
		2 \, \mu \int_\Omega \varepsilon(u) \dprod \varepsilon (v) \dx
		=
		a(u,v)
		,
		\\
		\dual{P_Q p}{q}
		&
		=
		\frac{1}{2 \, \mu} \int_\Omega p \, q \dx 
		=
		\frac{\lambda}{2 \, \mu} 
		c(p,q)
		.
	\end{aligned}
\end{equation}
Clearly, \eqref{eq:elasticity_preconditioner_blocks} satisfies our conditions of dimensional compatibility.
Essentially the same preconditioner but with $\dual{P_Q p}{q} = \lambda \, c(p,q)$, has been considered in \cite{Klawonn:1998:1}, where the emphasis was on robustness w.r.t.\ the parameter $\lambda$, which goes to infinity in the incompressible limit; see also \cite{KuchtaMardalMortensen:2019:2}.
In fact, it is straightforward to verify that with the scaling as in \eqref{eq:elasticity_preconditioner_blocks}, the constants in \eqref{eq:error_estimate_constants} are robust with respect to \emph{both} Lamé parameters $\mu$ and $\lambda$.
Indeed, using $\norm{\div u}_{L^2(\Omega)} = \norm{\trace \varepsilon(u)}_{L^2(\Omega)} \le \norm{\varepsilon(u)}_{L^2(\Omega;\R^{3 \times 3})}$, one easily finds
\begin{equation}
	\label{eq:elasticity_relevant_constants}
	\norm{A} = \alpha = 1,
	\quad 
	\norm{B} \le 1,
	\quad
	\beta > 0
\end{equation}
for all $(\mu,\lambda)$ which render the saddle-point system itself well-posed, i.e., $\mu > 0$ and $2 \, \mu + 3 \, \lambda > 0$; see for instance \cite[Proposition~3.13]{MarsdenHughes:1994:1}.
Notice that $\beta$ is the same inf-sup constant one obtains in the unscaled setting, i.e.\ with the norms $\norm{\varepsilon(u)}_{L^2(\Omega;\R^{3 \times 3})}$ and $\norm{p}_{L^2(\Omega)}$,
\begin{equation}
	\label{eq:elasticity_beta_unscaled}
	\beta
	=
	\inf_{p \in Q \setminus \{0\}} \sup_{u \in V \setminus \{0\}} \frac{\dual{Bu}{p}}{\norm{\varepsilon(u)}_{L^2(\Omega;\R^{3 \times 3})} \, \norm{p}_{L^2(\Omega)}}
\end{equation}
since 
\begin{equation*}
	\norm{u}_{P_V} \norm{p}_{P_Q}	
	=
	\sqrt{2 \, \mu} \, \norm{\varepsilon(u)}_{L^2(\Omega;\R^{3 \times 3})} \sqrt{1/(2 \, \mu)} \, \norm{p}_{P_Q}	
\end{equation*}
holds and thus the denominator in \eqref{eq:elasticity_beta_unscaled} can be replaced by $\norm{u}_{P_V} \norm{p}_{P_Q}$.
Consequently, the dimensional consistency of the preconditioner \eqref{eq:elasticity_preconditioner_blocks} also leads to a robust preconditioning of the saddle-point system \eqref{eq:elasticity_forms}.

We mention that the robust estimates in \eqref{eq:elasticity_relevant_constants} carry over to appropriate discrete settings.
In our experiments in \cref{subsec:experiments_elasticity}, we are using the Taylor--Hood finite element pair, \ie, we replace $V$ by a subspace $V_h$ spanned by piecewise trilinear functions on a geometrically conforming simplicial grid, and we replace $Q$ by a subspace $Q_h$ spanned by piecewise linear functions.
It is well known that this discretization is inf-sup stable on quasi-uniform mesh families, \ie, $\beta$ in \eqref{eq:elasticity_beta_unscaled} is bounded away from zero uniformly with respect to the mesh parameter; see for instance \cite{Verfuerth:1984:1,Wieners:2003:1} or \cite[Chapter~II, §4.2]{GiraultRaviart:1986:1}.

\subsection{Distributed Optimal Control of the Poisson Equation}
\label{subsec:distributed_optimal_control_Poisson}

In this section we consider a standard distributed optimal control problem for the Poisson equation:
\begin{equation}
	\label{eq:optimal_control_problem_Poisson}
	\begin{aligned}
		\text{Minimize} \quad & \frac{\beta}{2} \norm{u - u_d}_{L^2(\Omega)}^2 + \frac{\alpha}{2} \norm{f}_{L^2(\Omega)}^2
		\\
		\text{subject to} \quad &
		\left\{
			\begin{aligned}
				- \kappa \, \laplace u & = f & & \text{in } \Omega,
				\\
				u & = 0 & & \text{on } \Gamma = \partial\Omega.
			\end{aligned}
		\right.
	\end{aligned}
\end{equation}
This problem, with $\beta = \kappa = 1$ and an emphasis on robustness w.r.t.\ $\alpha$ was considered in \cite{SchoeberlZulehner:2007:1} and \cite[Section~4.1]{Zulehner:2011:1}.
Here we discuss a setting with physically meaningful constants and we show that the preconditioner developed in \cite[Section~4.1]{Zulehner:2011:1} can be extended to ensure robustness of w.r.t.\ \emph{all} problem parameters $\alpha,\beta,\kappa > 0$.
The corresponding scaling of the preconditioner building blocks is obtained from considerations of dimensional consistency.

The Poisson equation models many physical processes, and we consider it here in its role in stationary heat conduction on a bounded domain $\Omega \subset \R^3$.
The primal (state) variable~$u$ represents temperature, measured in $\unit{u} = \si{\kelvin}$.
The same is true for the desired state~$u_d$.
The heat conductivity~$\kappa$ with $\unit{\kappa} = \si{\watt\per\meter\per\kelvin}$ is a positive constant.
The control function~$f$ represents heating power, measured in $\unit{f} = \si{\watt\per\meter\cubed}$.
Finally, suppose that $\si{\obj}$ is the unit in which we measure the values of the objective.
It will become apparent that its choice does not matter.
Consequently, the coefficients $\beta > 0$ and $\alpha > 0$ have units $\unit{\beta} = \si{\obj\per\kelvin\squared\per\meter\cubed}$ and $\unit{\alpha} = \si{\obj\per\watt\squared\meter\cubed}$.

It is well known that \eqref{eq:optimal_control_problem_Poisson} has a unique solution $f \in L^2(\Omega)$ with associated weak solution of the state equation $u \in H^1_0(\Omega)$.
This solution is characterized via the KKT conditions associated with \eqref{eq:optimal_control_problem_Poisson}, which involves a unique adjoint state $p \in H^1_0(\Omega)$; see for instance \cite[Chapter~2.8]{Troeltzsch:2010:1}.
Since the control and adjoint state are related via $f = \alpha^{-1} p$, we can eliminate the control and obtain, similar to \cite[Section~4.1]{Zulehner:2011:1}, a saddle-point system \eqref{eq:self-adjoint_saddle-point_system} with
\begin{equation}
	\label{eq:optimal_control_problem_Poisson_forms}
	\begin{aligned}
		a(u, v) 
		& 
		= 
		\beta \int_\Omega u \, v \dx, 
		&
		b(u, q) 
		& 
		= 
		\kappa \int_\Omega \nabla u \cdot \nabla q \dx, 
		\\
		c(p,q) 
		&
		= 
		\alpha^{-1} \int_\Omega p \, q \dx,
		&
		f(v) 
		&
		= 
		\beta \int_\Omega u_d \, v \dx
	\end{aligned}
\end{equation}
as well as $g = 0$, where $u, v \in V = H^1_0(\Omega)$ and $p, q \in Q = H^1_0(\Omega)$.
Since the unit of the Lagrangian \eqref{eq:Lagrangian} is $\unit{\cL(u,p)} = \si{\obj}$, the physical units associated with the spaces $V$ and $Q$ and their duals are
\begin{equation*}
	\begin{aligned}
		\unit{V} & = \si{\kelvin},
						 &
		\unit{V^*} & = \si{\obj\per\kelvin},
		\\
		\unit{Q} & = \si{\obj\per\watt},
						 &
		\unit{Q^*} & = \si{\watt}.
	\end{aligned}
\end{equation*}

Following \cite[Section~3]{Zulehner:2011:1}, we continue the discussion replacing $V$ and $Q$ by finite-dimensional subspaces $V_h$ and $Q_h$ in order to avoid technicalities.
Suppose that $\tA$, $\tB$ and $\tC$ are matrices representing the operators $A$, $B$ and $C$, respectively, with respect to the chosen bases of the subspaces  $V_h$ and $Q_h$.
Clearly, $\tA$ and $\tC$ are symmetric and positive definite, and so are their negative Schur complements $\tB \tA^{-1} \tB^\transp$ and $\tB^\transp \tC^{-1} \tB$.
Therefore, the theory in \cite[Section~3.3]{Zulehner:2011:1} applies, which reveals that
\begin{equation}
	\label{eq:optimal_control_problem_preconditioner_family}
	\tP_V = \tA + \interp[big]{\tA}{\tB^\transp \tC^{-1} \tB}{\theta}, 
	\quad
	\tP_Q = \tC + \interp[big]{\tC}{\tB \tA^{-1} \tB^\transp}{1-\theta} 
\end{equation}
is a class of robust preconditioners for all $\theta \in [0,1]$.
Here 
\begin{equation*}
	\interp{\tV}{\tW}{\theta} \coloneqq \tV^{1/2} \paren[big](){\tV^{-1/2} \tW \tV^{-1/2}}^\theta \tV^{1/2} 
\end{equation*}
denotes the interpolation between symmetric, positive definite matrices~$\tV$ and $\tW$ of equal size.

Notice that although the emphasis in \cite{Zulehner:2011:1} was on showing the $\alpha$-robustness (in addition to robustness w.r.t.\ the discretization) for this preconditioner in case $\beta = \kappa = 1$, in fact it is, by construction, robust w.r.t.\ \emph{all} problem parameters $\alpha,\beta,\kappa > 0$.

Let us show that \eqref{eq:optimal_control_problem_preconditioner_family} is meaningful from the viewpoint of physical dimensions.
To this end, we note that $\tA$, $\tB$ and $\tC$ map between coefficient vectors w.r.t.\ the chosen bases.
In order for these matrices to retain the same mapping properties of their continuous counterparts $A$, $B$ and $C$, respectively, we leave the dimensions to the coefficient vectors and choose a dimensionless basis $\{\varphi_j\}$ of $V_h$ and a dimensionless basis $\{\psi_j\}$ of $Q_h$, i.e., $\unit{\varphi_j} = \unit{\psi_j} = 1$.
Then $\tA_{ij} \coloneqq \dual{A \varphi_j}{\varphi_i}$, $\tB_{ij} \coloneqq \dual{B \varphi_j}{\psi_i}$ and $\tC_{ij} \coloneqq \dual{C \psi_j}{\psi_i}$ hold.
Therefore, we obtain $\unit{\tA} = \frac{\unit{V^*}}{\unit{V}} = \si{\obj\per\kelvin\squared}$, $\unit{\tB} = \frac{\unit{Q^*}}{\unit{V}} = \si{\watt\per\kelvin}$, $\unit{\tB^\transp} = \frac{\unit{V^*}}{\unit{Q}} = \si{\watt\per\kelvin}$ and $\unit{\tC} = \frac{\unit{Q^*}}{\unit{Q}} = \si{\watt\squared\per\obj}$.
This implies that $\unit{\tA} = \unit{\tB^\transp \tC^{-1} \tB}$ and $\unit{\tC} = \unit{\tB \tA^{-1} \tB^\transp}$, so that \eqref{eq:optimal_control_problem_preconditioner_family} is meaningful from a dimensional point of view for all $\theta \in [0,1]$.

The preconditioner \eqref{eq:optimal_control_problem_preconditioner_family} may seem difficult to implement in practice.
However, as was pointed out in \cite[Section~4.1]{Zulehner:2011:1}, the case $\theta = 1/2$ leads to a simple representation in case that $V_h$ and $Q_h$ are the same space and the same basis $\varphi_j = \psi_j$ is chosen.
(A standard example is to consider a space of continuous, piecewise linear finite element functions.)
In this case, one has 
\begin{equation*}
	\tA = \beta \, \tM,
	\quad
	\tB = \kappa \, \tK,
	\quad
	\tC = \alpha^{-1} \, \tM,
\end{equation*}
where 
\begin{equation*}
	\tM_{ij} 
	=
	\int_\Omega \varphi_j \, \varphi_i \dx,
	\quad
	\tK_{ij}
	=
	\int_\Omega \nabla \varphi_j \cdot \nabla \psi_i \dx.
\end{equation*}
Notice that $\unit{\tM} = \si{\meter\cubed}$ and $\unit{\tK} = \si{\meter}$ and that, clearly, both matrices are symmetric.
Moreover, similarly to \cite[Section~4.1]{Zulehner:2011:1} we observe that
\begin{equation*}
	\interp[big]{\tM}{\tK^\transp \tM^{-1} \tK}{1/2}
	=
	\tM^{1/2} \paren[Big](){\underbrace{\tM^{-1/2} \tK^\transp \tM^{-1} \tK \tM^{-1/2}}_{= \paren[auto](){\tM^{-1/2} \tK \tM^{-1/2}}^2}}^{1/2} \tM^{1/2} 
	=
	\tK
\end{equation*}
and likewise,
\begin{equation*}
	\interp[big]{\tM}{\tK \tM^{-1} \tK^\transp}{1/2}
	=
	\tK
	.
\end{equation*}
Taking into account that $\interp{\gamma \tV}{\delta \tW}{\theta} = \gamma^{1-\theta} \delta^\theta [\tV, \tW]_\theta$ holds for $\gamma, \delta > 0$, we find that the preconditioner in \eqref{eq:optimal_control_problem_preconditioner_family} indeed has the following simple representation in case $\theta = 1/2$, 
\begin{equation}
	\label{eq:optimal_control_problem_Poisson_preconditioner_particular}
	\tP_V = \beta \, \tM + (\alpha \, \beta)^{1/2} \kappa \, \tK,
	\quad
	\tP_Q = \alpha^{-1} \tM + (\alpha \, \beta)^{-1/2} \kappa \, \tK
	=
	(\alpha \, \beta)^{-1} \, \tP_V
	.
\end{equation}
As was mentioned before, this preconditioner is dimensionally consistent and robust w.r.t.\ \emph{all} problem parameters $\alpha,\beta,\kappa > 0$.

\subsection{Distributed Optimal Control of the Stokes Equation}
\label{subsec:distributed_optimal_control_Stokes}

In this section we apply the general framework of \cref{sec:general_framework} to a nested saddle-point problem.
Compared with the example in \cref{subsec:distributed_optimal_control_Poisson}, the state equation is no longer elliptic but has a saddle-point structure in its own right.
Our example is a distributed optimal control example of the Stokes system similar to \cite[Section~4.2]{Zulehner:2011:1} on a bounded domain $\Omega \subset \R^3$, but equipped with a full set of constants rendering the problem description physically meaningful:
\begin{equation}
\label{eq:optimal_control_problem_Stokes}
	\begin{aligned}
		\text{Minimize} \quad & \frac{\beta}{2} \norm{u - u_d}_{L^2(\Omega)^3}^2 + \frac{\alpha}{2} \norm{f}_{L^2(\Omega)^3}^2
		\\
		\text{subject to} \quad &
		\left\{
			\begin{aligned}
				- \mu \, \laplace u + \nabla p & = f & & \text{in } \Omega,
				\\
				\div u & = 0 & & \text{in } \Omega,
				\\
				u & = 0 & & \text{on } \Gamma = \partial\Omega.
			\end{aligned}
		\right.
	\end{aligned}
\end{equation}
In \cite{Zulehner:2011:1}, this problem was considered with $\beta = \mu = 1$.
As in \eqref{eq:Stokes_system_with_bilinear_forms}, the unit of the dynamic viscosity is $\unit{\mu} = \si{\newton\second\per\meter\squared}$ while $\unit{u} = \si{\meter\per\second}$ and $\unit{p} = \si{\newton\per\meter\squared}$ as well as $\unit{f} = \si{\newton\per\meter\cubed}$.
Moreover, the coefficients $\beta > 0$ and $\alpha > 0$ have units $\unit{\beta} = \si{\obj\per\meter\tothe{5}\second\squared}$ and $\unit{\alpha} = \si{\obj\per\newton\squared\meter\cubed}$.
The desired velocity state finally has $\unit{u_d} = \si{\meter\per\second}$.

Similar to \cref{subsec:distributed_optimal_control_Poisson}, we formulate the KKT optimality conditions characterizing the unique optimal control $f \in L^2(\Omega)^3$ with associated weak solution $(u,p) \in H^1_0(\Omega)^3 \times L^2_0(\Omega)$ of the Stokes system.
The optimality system involves an adjoint state $(w,r) \in H^1_0(\Omega)^3 \times L^2_0(\Omega)$.
Recall that $L^2_0(\Omega)$ denotes the space of $L^2(\Omega)$ with zero mean.
An elimination of the optimal control via the relation $f = \alpha^{-1} w$ leads to a saddle-point system \eqref{eq:self-adjoint_saddle-point_system} with $g = 0$ and 
\makeatletter
\begin{equation}
	\label{eq:optimal_control_problem_Stokes_forms}
	\begin{aligned}
		a((u,p), (v,q)) 
		& 
		= 
		\beta \int_\Omega u \cdot v \dx, 
		\ltx@ifpackageloaded{dgruyter-numapde}{%
			\\
			}{%
			&
		}
		b((u,p), (z,s)) 
		& 
		= 
		\mu \int_\Omega \nabla u \cdot \nabla z \dx - \int_\Omega s \, \div u \dx - \int_\Omega p \, \div z \dx, 
		\\
		c((w,r),(z,s)) 
		&
		= 
		\alpha^{-1} \int_\Omega w \cdot z \dx,
		\ltx@ifpackageloaded{dgruyter-numapde}{%
			\\
			}{%
			&
		}
		f((v,q)) 
		&
		= 
		\beta \int_\Omega u_d \cdot v \dx,
	\end{aligned}
\end{equation}
\makeatother
compare \cite[Section~4.2]{Zulehner:2011:1}.
The relevant spaces are $V = Q = H^1_0(\Omega)^3 \times L^2_0(\Omega)$, equipped with the units
\begin{equation*}
	\begin{aligned}
		\unit{V} & = \begin{pmatrix} \si{\meter\per\second} \\ \si{\newton\per\meter\squared} \end{pmatrix},
						 &
		\unit{V^*} & = \begin{pmatrix} \si{\obj\second\per\meter} \\ \si{\obj\meter\squared\per\newton} \end{pmatrix},
		\\
		\unit{Q} & = \begin{pmatrix} \si{\obj\per\newton} \\ \si{\obj\second\per\meter\cubed} \end{pmatrix},
						 &
		\unit{Q^*} & = \begin{pmatrix} \si{\newton} \\ \si{\meter\cubed\per\second} \end{pmatrix}.
	\end{aligned}
\end{equation*}
We discretize the spaces $V$ and $Q$ by identical Taylor--Hood finite element spaces.
This then leads to a discretized optimality system of the form \eqref{eq:self-adjoint_saddle-point_system} with
\begin{equation*}
	\begin{aligned}
		\tA
		=
		\begin{bmatrix}
			\beta \, \tM & \\ & 0
		\end{bmatrix},
		\quad
		\tB
		=
		\begin{bmatrix}
			\mu \, \tK & -\tD^\transp \\ -\tD &
		\end{bmatrix},
		\quad
		\tC
		=
		\begin{bmatrix}
			\alpha^{-1} \tM & \\ & 0
		\end{bmatrix}
		=
		(\alpha \, \beta)^{-1} \tA
		.
	\end{aligned}
\end{equation*}
Here $\tM$ and $\tK$ are the mass and stiffness matrices over the vector-valued, quadratic finite element space, and $\tD$ is the matrix representation of the bilinear form $\int_\Omega r \, \div u \dx$.
As in the previous section we have $\unit{\tM} = \si{\meter\cubed}$ and $\unit{\tK} = \si{\meter}$ and, moreover, $\unit{\tD} = \si{\meter\squared}$.
Repeating and slightly extending the analysis in \cite[Section~4.2]{Zulehner:2011:1} by working in the problem parameters $\beta$ and $\mu$ shows that
\makeatletter
\begin{equation}
	\label{eq:optimal_control_problem_Stokes_preconditioner_particular}
	\begin{aligned}
		\tP_V 
		&
		= 
		\begin{bmatrix}
			\beta \, \tM + (\alpha \, \beta)^{1/2} \mu \, \tK & \\ & \alpha \, \beta \, \tD \paren[big][]{\beta \, \tM + (\alpha \, \beta)^{1/2} \mu \, \tK}^{-1} \tD^\transp
		\end{bmatrix},
		\ltx@ifpackageloaded{dgruyter-numapde}{%
			\\
			}{%
			&
		}
		\tP_Q 
		&
		= 
		(\alpha \, \beta)^{-1} \tP_V
	\end{aligned}
\end{equation}
\makeatother
is a preconditioner for the problem at hand, which is robust not only w.r.t.\ $\alpha$ but \emph{all} problem parameters $\alpha,\beta,\mu > 0$.
Once again, the preconditioner is also dimensionally consistent, as is easily checked.

\section{Numerical Experiments}
\label{sec:numerical_experiments}

In this section we describe a number of experiments with the purpose of verifying numerically the robustness and mesh independence of the preconditioners for model problems considered in \cref{subsect:introductory_example,subsec:elasticity,subsec:distributed_optimal_control_Poisson,subsec:distributed_optimal_control_Stokes}.
In each case, we report the convergence history of preconditioned \minres for the norm of the residual,
\begin{equation}
	\label{eq:total_residual_norm}
	\norm{r}_{P^{-1}}
	\coloneqq
	\paren[big](){\norm{r_1}_{P_V^{-1}}^2 + \norm{r_2}_{P_Q^{-1}}^2}^{1/2}
	.
\end{equation}
In each case, we stop the iteration once this quantity has been reduced by a factor of $10^{-6}$ compared to its initial value, which is the residual associated with all-zero initial guesses.
We use the \minres implementation in \matlab described in \cite{HerzogSoodhalter:2017:1} and available from \cite{HerzogSoodhalter:2016:2}, which allows us to monitor the convergence histories of $\norm{r_1}_{P_V^{-1}}$ and $\norm{r_2}_{P_Q^{-1}}$ separately.

For all problems, we assembled the matrices and right hand side vectors using the finite element package \fenics (version~2019.1); see \cite{LoggMardalWells:2012:1}.
In each case, a geometrically conforming tetrahedral grid was generated for the coarse mesh level and then refined for the finer levels. 
The matrices and vectors were then exported and read into \matlab (version~R2020a).

All preconditioners discussed in this paper, i.e., \eqref{eq:Stokes_preconditioner_blocks}, \eqref{eq:elasticity_preconditioner_blocks}, \eqref{eq:optimal_control_problem_Poisson_preconditioner_particular} and \eqref{eq:optimal_control_problem_Stokes_preconditioner_particular}, were discussed in their ideal forms only.
This is since our main emphasis was on their robustness w.r.t.\ various problem parameters, in line with their dimensional consistency.
In our experiments we indeed employ these preconditioners in their ideal forms and apply them using a Cholesky factorization with permutations to reduce fill-in, as provided by \matlab's \texttt{chol} command.
For truly large-scale problems one would replace the preconditioner's building blocks with spectrally equivalent implementations, e.g., based on geometric multigrid approaches.
This is well known and it does not interfere with our considerations of dimensional consistency.

\subsection{Stokes}
\label{subsec:experiments_Stokes}

Our Stokes model problem is a 3D variant of \cite[Example~3.1.1]{ElmanSilvesterWathen:2014:1}.
We have $\Omega = (-1,1)^3 \si{\meter\cubed}$ and the boundary $\Gamma$ is split into three parts,
\begin{equation*}
	\begin{aligned}
		\Gamma_\text{inflow} & = \setDef{x \in \Gamma}{x_1 = \SI{-1}{\meter}},
		\\
		\Gamma_\text{noslip} & = \setDef{x \in \Gamma}{\abs{x_2} = \SI{1}{\meter}} \cup \setDef{x \in \R^3}{\abs{x_3} = \SI{1}{\meter}},
		\\
		\Gamma_\text{outflow} & = \setDef{x \in \Gamma}{x_1 = \SI{1}{\meter}}.
	\end{aligned}
\end{equation*}
We impose the conditions $u(x) = \paren[auto](){(1 - x_2^2) \, (1 - x_3^2), 0, 0} \si{\meter\per\second}$ on $\Gamma_\text{inflow}$, $u(x) = (0, 0, 0)~\si{\meter\per\second}$ on $\Gamma_\text{noslip}$ and homogeneous natural (do-nothing) outflow boundary conditions on $\Gamma_\text{outflow}$.
Moreover, we use a volume force density of $F = (0,0,0)~\si{\newton\per\meter\cubed}$ in \eqref{eq:Stokes_forms}.
The problem was discretized with the Taylor--Hood finite element pair, i.e., we used continuous, piecewise quadratic functions for the components of the velocity~$u$ and continuous, piecewise linear functions for the pressure~$p$.

\begin{figure}[htp]
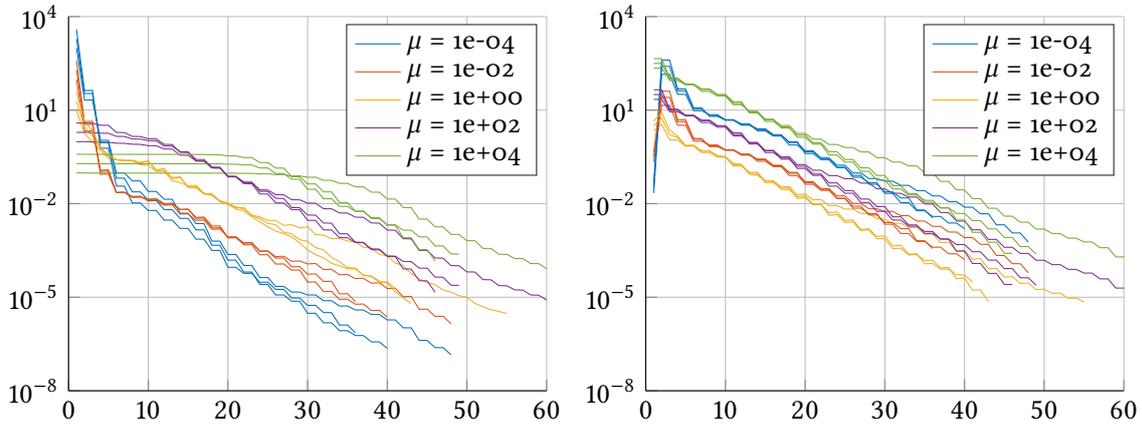

	\begin{subfigure}[b]{0.42\textwidth}
		\input{experiments/residuals_1_Stokes3D.tex}
	\end{subfigure}
	\hspace{0.05\textwidth}
	\begin{subfigure}[b]{0.42\textwidth}
		\input{experiments/residuals_2_Stokes3D.tex}
	\end{subfigure}
	\caption{Convergence of the norms $\norm{r_1}_{P_V^{-1}}$ (left) and $\norm{r_2}_{P_Q^{-1}}$ (right) of the residuals $r_1 \in V^*$ and $r_2 \in Q^*$ for the three discretizations of the Stokes problem (\cref{subsec:experiments_Stokes}) with preconditioner \eqref{eq:Stokes_preconditioner_blocks}.}
	\label{fig:residuals_Stokes3D}
\end{figure}

Although this model problem uses slightly more general boundary conditions than \eqref{eq:Stokes_system_with_bilinear_forms}, this does not affect the dimensional consistency of the preconditioner \eqref{eq:Stokes_preconditioner_blocks} nor the proof of robustness w.r.t.\ the viscosity~$\mu$ and the mesh size, based on \eqref{eq:error_estimate_constants}.
The convergence results reported in \cref{fig:residuals_Stokes3D,tab:convergence_Stokes3D} confirm this.

\begin{table}[htpb]
	\centering
	\begin{tabular}{@{}rrrrrrrrr@{}}
  \toprule
           level &       $\dim(V)$ &       $\dim(Q)$ & $\mu = 10^{-4}$ & $\mu = 10^{-2}$ &  $\mu = 10^{0}$ &  $\mu = 10^{2}$ &  $\mu = 10^{4}$ & \si{\newton\second\per\meter\squared} \\
  \midrule
  \texttt{    1} & \texttt{  2187} & \texttt{   125} & \texttt{    48} & \texttt{    48} & \texttt{    55} & \texttt{    60} & \texttt{    60} &                                       \\
  \texttt{    2} & \texttt{ 14739} & \texttt{   729} & \texttt{    40} & \texttt{    40} & \texttt{    43} & \texttt{    46} & \texttt{    46} &                                       \\
  \texttt{    3} & \texttt{107811} & \texttt{  4913} & \texttt{    36} & \texttt{    36} & \texttt{    41} & \texttt{    49} & \texttt{    49} &                                       \\
  \bottomrule
\end{tabular}

	\caption{Iteration numbers for the Stokes problem (\cref{subsec:experiments_Stokes}) required to reach a relative reduction by $10^{-6}$ of the initial residual norm \eqref{eq:total_residual_norm}.}
	\label{tab:convergence_Stokes3D}
\end{table}

\subsection{Nearly Incompressible Elasticity}
\label{subsec:experiments_elasticity}

We consider a rod of square cross section $\Omega = (0,100) \times (0,10) \times (0,10)~\si{\milli\meter\cubed}$, whose boundary~$\Gamma$ is decomposed into
\begin{equation*}
	\begin{aligned}
		\Gamma_D & = \setDef{x \in \Gamma}{x_1 = \SI{0}{\milli\meter}},
		\\
		\Gamma_N & = \Gamma \setminus \Gamma_D.
	\end{aligned}
\end{equation*}
Homogeneous Dirichlet conditions for the displacement $u$ are imposed at the clamping boundary~$\Gamma_D$.
On~$\Gamma_N$ we have natural (traction) boundary conditions.
The imposed traction pressure is $F = (0,0,0)~\si{\newton\per\milli\meter\squared}$ except at the forcing boundary located at $x_1 = \SI{100}{\milli\meter}$, where a uniform pressure of $F = (1,0,0)~\si{\newton\per\milli\meter\squared}$ is imposed.
As in the Stokes problem from \cref{subsec:experiments_Stokes}, we used continuous, piecewise quadratic functions for the displacement and continuous, piecewise linear functions for the pressure.

\begin{figure}[htp]
	\begin{subfigure}[b]{0.42\textwidth}
%
\definecolor{mycolor1}{rgb}{0.00000,0.44700,0.74100}%
\definecolor{mycolor2}{rgb}{0.85000,0.32500,0.09800}%
\definecolor{mycolor3}{rgb}{0.92900,0.69400,0.12500}%
\definecolor{mycolor4}{rgb}{0.49400,0.18400,0.55600}%
\definecolor{mycolor5}{rgb}{0.46600,0.67400,0.18800}%
\begin{tikzpicture}

\begin{axis}[%
width=0.951\linewidth,
height=0.75\linewidth,
at={(0\linewidth,0\linewidth)},
scale only axis,
xmin=0,
xmax=30,
ymode=log,
ymin=1e-08,
ymax=100000,
yminorticks=true,
axis background/.style={fill=white},
axis x line*=bottom,
axis y line*=left,
xmajorgrids,
ymajorgrids,
yminorgrids,
legend style={legend cell align=left, align=left, draw=white!15!black}
]
\addplot [color=mycolor1]
  table[row sep=crcr]{%
1	9447.95028583424\\
2	4722.8395435662\\
3	4624.73256862197\\
4	59.0937405080328\\
5	59.2465937922135\\
6	3.67341667814458\\
7	3.67386005603341\\
8	0.186376619037188\\
9	0.186378069415916\\
10	0.0210602125377004\\
11	0.0210602324865092\\
12	0.00315764233787065\\
13	0.0031576429487388\\
14	0.00058515629707197\\
15	0.000585156325398421\\
16	0.000103927328573964\\
17	0.000103927332298738\\
18	1.74963687185437e-05\\
19	1.7496370005952e-05\\
20	1.60486230376741e-06\\
21	1.60486339896012e-06\\
22	2.62862454602611e-07\\
23	2.62863717000255e-07\\
24	5.93985119993403e-08\\
25	5.94002623122695e-08\\
26	1.44121928088073e-08\\
27	1.44133089998747e-08\\
28	2.70033705958234e-09\\
};
\addlegendentry{$\mu\text{ = 1e-04, }\lambda\text{ = 1e+04}$}

\addplot [color=mycolor2]
  table[row sep=crcr]{%
1	944.795028583424\\
2	472.283954356567\\
3	462.46692379362\\
4	5.90897915018874\\
5	5.92506191007471\\
6	0.367273988116382\\
7	0.367363033115989\\
8	0.0186326371105679\\
9	0.0186449635133425\\
10	0.00210633149297392\\
11	0.00211701088131224\\
12	0.000317330167275963\\
13	0.000330154582021531\\
14	6.11630662301548e-05\\
15	7.1901508765193e-05\\
16	1.27660485972807e-05\\
17	2.09196018086566e-05\\
18	3.52064280183401e-06\\
19	6.97263370270215e-06\\
20	6.39371216769687e-07\\
21	1.9605479886269e-06\\
22	3.21032817842458e-07\\
23	8.73473550492927e-07\\
24	1.97234668524879e-07\\
25	4.96209366410586e-07\\
26	1.19407814561282e-07\\
27	2.15278309118708e-07\\
28	2.9349687982054e-08\\
};
\addlegendentry{$\mu\text{ = 1e-02, }\lambda\text{ = 1e+02}$}

\addplot [color=mycolor3]
  table[row sep=crcr]{%
1	94.4795028583424\\
2	47.2283954356619\\
3	33.8994509225183\\
4	0.405439857499053\\
5	1.78976874745568\\
6	0.0760771643615527\\
7	0.206236217060657\\
8	0.00778545580312549\\
9	0.0223424973778946\\
10	0.00177356278140661\\
11	0.00340172897077071\\
12	0.00034880996099828\\
13	0.000559430362396652\\
14	6.50913849752935e-05\\
15	9.84177751471921e-05\\
16	1.1083202331382e-05\\
};
\addlegendentry{$\mu\text{ = 1e+00, }\lambda\text{ = 1e+00}$}

\addplot [color=mycolor4]
  table[row sep=crcr]{%
1	9.44795028583424\\
2	4.72283954356581\\
3	0.000661526312454681\\
4	1.49106666257572e-05\\
5	6.10935991918175e-09\\
};
\addlegendentry{$\mu\text{ = 1e+02, }\lambda\text{ = 1e-02}$}

\addplot [color=mycolor5]
  table[row sep=crcr]{%
1	0.944795028583424\\
2	0.472283954356627\\
3	6.6155937588943e-09\\
};
\addlegendentry{$\mu\text{ = 1e+04, }\lambda\text{ = 1e-04}$}

\addplot [color=mycolor1, forget plot]
  table[row sep=crcr]{%
1	8258.7108671638\\
2	4129.29778050244\\
3	4043.60394546332\\
4	43.0877065993229\\
5	43.1792682871825\\
6	3.73920500564319\\
7	3.73983571687524\\
8	0.346984065185103\\
9	0.346991139290355\\
10	0.0732409472721828\\
11	0.0732412843827095\\
12	0.0103745691920289\\
13	0.0103745750144631\\
14	0.0013735919660138\\
15	0.001373592105264\\
16	0.000269761971212037\\
17	0.000269761980435564\\
18	3.49205375412106e-05\\
19	3.49205394184846e-05\\
20	1.13680710842019e-05\\
21	1.13680732172254e-05\\
22	2.44931940324064e-06\\
23	2.44932038237127e-06\\
24	3.18257030611842e-07\\
25	3.18258395032997e-07\\
26	1.2511376650686e-07\\
27	1.25115719349113e-07\\
28	2.52431791319269e-08\\
29	2.52440378417082e-08\\
30	2.36210828478591e-09\\
};
\addplot [color=mycolor2, forget plot]
  table[row sep=crcr]{%
1	825.87108671638\\
2	412.929778050264\\
3	404.354869587947\\
4	4.30848273223623\\
5	4.31821422467136\\
6	0.373854894922931\\
7	0.373971289480515\\
8	0.0346889278835397\\
9	0.0347056604031481\\
10	0.00732351419221207\\
11	0.00733602958239995\\
12	0.00103877919693599\\
13	0.00104834970640361\\
14	0.000138760311104682\\
15	0.000149282584281342\\
16	2.93090437923148e-05\\
17	3.75136892240268e-05\\
18	4.85453407118213e-06\\
19	1.02926308653727e-05\\
20	3.34976848337381e-06\\
21	7.2344511888216e-06\\
22	1.55788601786319e-06\\
23	2.65308326889915e-06\\
24	3.44604051211009e-07\\
25	9.90042258667396e-07\\
26	3.89082988540764e-07\\
27	7.98348953510353e-07\\
28	1.60800311203844e-07\\
29	2.62688235896583e-07\\
30	2.10761697788392e-08\\
};
\addplot [color=mycolor3, forget plot]
  table[row sep=crcr]{%
1	82.587108671638\\
2	41.2929778050247\\
3	29.6560863194073\\
4	0.295784960732549\\
5	1.43480548457655\\
6	0.0863688658834231\\
7	0.192004927662989\\
8	0.0125392532982465\\
9	0.0256550983031618\\
10	0.00360758538342272\\
11	0.00497755438501742\\
12	0.000428654249264779\\
13	0.00076643259435437\\
14	6.69284120653895e-05\\
15	0.000116567739260391\\
16	1.47232539659037e-05\\
};
\addplot [color=mycolor4, forget plot]
  table[row sep=crcr]{%
1	8.2587108671638\\
2	4.12929778050189\\
3	0.00057908388674532\\
4	1.09156686765623e-05\\
5	4.91968344920501e-09\\
};
\addplot [color=mycolor5, forget plot]
  table[row sep=crcr]{%
1	0.82587108671638\\
2	0.412929778050245\\
3	5.7911285077092e-09\\
};
\addplot [color=mycolor1, forget plot]
  table[row sep=crcr]{%
1	7663.60742996663\\
2	3831.80075930102\\
3	3752.39403778256\\
4	37.5252069055682\\
5	37.6003571098045\\
6	3.64508875600269\\
7	3.64577703551259\\
8	0.347830584138542\\
9	0.347837362448994\\
10	0.0575532753391978\\
11	0.057553500043811\\
12	0.00875235798788731\\
13	0.0087523622390063\\
14	0.000734595059932713\\
15	0.000734595096794258\\
16	0.000102742348252553\\
17	0.00010274235155364\\
18	1.46856966009973e-05\\
19	1.46856978321769e-05\\
20	3.46969978606526e-06\\
21	3.46970125956153e-06\\
22	8.19648609242833e-07\\
23	8.19649549486162e-07\\
24	1.03609154655002e-07\\
25	1.03610359005111e-07\\
26	3.92391750154026e-08\\
27	3.92411920150773e-08\\
28	1.08038038480661e-08\\
29	1.08047401253589e-08\\
30	2.30055263064748e-09\\
};
\addplot [color=mycolor2, forget plot]
  table[row sep=crcr]{%
1	766.360742996663\\
2	383.180075930113\\
3	375.234283731515\\
4	3.75226996841431\\
5	3.76028897565023\\
6	0.364445831801504\\
7	0.364569661500017\\
8	0.0347733645597286\\
9	0.0347875670388904\\
10	0.00575449221266459\\
11	0.00576593896750341\\
12	0.000876558419296019\\
13	0.000884980160486146\\
14	7.4255043668905e-05\\
15	8.15560807456547e-05\\
16	1.14039326455507e-05\\
17	1.82588572239095e-05\\
18	2.60907075102885e-06\\
19	5.9733978366181e-06\\
20	1.41090262240911e-06\\
21	3.41263174307456e-06\\
22	8.05837657736301e-07\\
23	1.47182785396437e-06\\
24	1.85958082891186e-07\\
25	5.3191610990737e-07\\
26	2.01235096036516e-07\\
27	4.45141405926546e-07\\
28	1.21085049136252e-07\\
29	1.8656400131826e-07\\
30	2.71152091787841e-08\\
};
\addplot [color=mycolor3, forget plot]
  table[row sep=crcr]{%
1	76.6360742996663\\
2	38.3180075929927\\
3	27.5270881309128\\
4	0.257663132723686\\
5	1.28896214757444\\
6	0.0870118728365883\\
7	0.180713167055412\\
8	0.0118363970116432\\
9	0.0247099846524509\\
10	0.00279728629947642\\
11	0.00429579901394065\\
12	0.00040502400278506\\
13	0.000701228415713192\\
14	3.87416392107026e-05\\
15	8.97782295800648e-05\\
16	8.88511732776663e-06\\
};
\addplot [color=mycolor4, forget plot]
  table[row sep=crcr]{%
1	7.66360742996663\\
2	3.83180075929947\\
3	0.000537629344444734\\
4	9.52060435565225e-06\\
5	4.41542971765229e-09\\
};
\addplot [color=mycolor5, forget plot]
  table[row sep=crcr]{%
1	0.766360742996663\\
2	0.38318007593012\\
3	5.37656287790364e-09\\
};
\end{axis}
\end{tikzpicture}%
	\end{subfigure}
	\hspace{0.05\textwidth}
	\begin{subfigure}[b]{0.42\textwidth}
%
\definecolor{mycolor1}{rgb}{0.00000,0.44700,0.74100}%
\definecolor{mycolor2}{rgb}{0.85000,0.32500,0.09800}%
\definecolor{mycolor3}{rgb}{0.92900,0.69400,0.12500}%
\definecolor{mycolor4}{rgb}{0.49400,0.18400,0.55600}%
\definecolor{mycolor5}{rgb}{0.46600,0.67400,0.18800}%
\begin{tikzpicture}

\begin{axis}[%
width=0.951\linewidth,
height=0.75\linewidth,
at={(0\linewidth,0\linewidth)},
scale only axis,
xmin=0,
xmax=30,
ymode=log,
ymin=1e-08,
ymax=100000,
yminorticks=true,
axis background/.style={fill=white},
axis x line*=bottom,
axis y line*=left,
xmajorgrids,
ymajorgrids,
yminorgrids,
legend style={legend cell align=left, align=left, draw=white!15!black}
]
\addplot [color=mycolor1]
  table[row sep=crcr]{%
1	0\\
2	4723.9750064234\\
3	3765.20545342184\\
4	671.526709756286\\
5	671.499561735944\\
6	167.814300874971\\
7	167.81428145361\\
8	37.8061192059432\\
9	37.8061191916426\\
10	12.708688961922\\
11	12.7086889618559\\
12	4.92097485064682\\
13	4.92097485064604\\
14	2.11838869784117\\
15	2.11838869784115\\
16	0.892759747240159\\
17	0.892759747240158\\
18	0.366305729032492\\
19	0.366305729032491\\
20	0.110940195374497\\
21	0.110940195374497\\
22	0.0448984492266969\\
23	0.0448984492266969\\
24	0.0213383956234456\\
25	0.0213383956234455\\
26	0.0104697160425109\\
27	0.0104697160425109\\
28	0.00386645998875398\\
};
\addlegendentry{$\mu\text{ = 1e-04, }\lambda\text{ = 1e+04}$}

\addplot [color=mycolor2]
  table[row sep=crcr]{%
1	0\\
2	472.39750064234\\
3	376.489965950325\\
4	67.1481598515382\\
5	67.1453032191665\\
6	16.7781691336391\\
7	16.7781652329078\\
8	3.77948566692191\\
9	3.7794855453419\\
10	1.27033998463324\\
11	1.2703399491287\\
12	0.491829533411381\\
13	0.491829516528262\\
14	0.211690139529955\\
15	0.211690132779965\\
16	0.089198978776656\\
17	0.0891989756974979\\
18	0.036592674879441\\
19	0.0365926738895517\\
20	0.0110808295551504\\
21	0.0110808292451598\\
22	0.00448392389732464\\
23	0.00448392375015578\\
24	0.00213071287638772\\
25	0.00213071277908591\\
26	0.00104522232498791\\
27	0.00104522229428964\\
28	0.000385929306519392\\
};
\addlegendentry{$\mu\text{ = 1e-02, }\lambda\text{ = 1e+02}$}

\addplot [color=mycolor3]
  table[row sep=crcr]{%
1	0\\
2	47.239750064234\\
3	14.3445300648091\\
4	4.00509396658015\\
5	1.07991736693233\\
6	0.424199275384055\\
7	0.145482598263676\\
8	0.0484150738889143\\
9	0.0139930254823847\\
10	0.00721272578887299\\
11	0.00207673160292578\\
12	0.0012276482822126\\
13	0.000306856215380093\\
14	0.000207685074621321\\
15	5.64259379792504e-05\\
16	3.64210980451618e-05\\
};
\addlegendentry{$\mu\text{ = 1e+00, }\lambda\text{ = 1e+00}$}

\addplot [color=mycolor4]
  table[row sep=crcr]{%
1	0\\
2	4.7239750064234\\
3	4.63298652051324e-08\\
4	9.81909892986719e-05\\
5	3.71551365616902e-13\\
};
\addlegendentry{$\mu\text{ = 1e+02, }\lambda\text{ = 1e-02}$}

\addplot [color=mycolor5]
  table[row sep=crcr]{%
1	0\\
2	0.47239750064234\\
3	3.9579236370769e-16\\
};
\addlegendentry{$\mu\text{ = 1e+04, }\lambda\text{ = 1e-04}$}

\addplot [color=mycolor1, forget plot]
  table[row sep=crcr]{%
1	0\\
2	4129.35543317943\\
3	3292.41828204048\\
4	536.546168823102\\
5	536.531351383396\\
6	158.353601514464\\
7	158.353571709348\\
8	48.2465449471761\\
9	48.2465448454173\\
10	22.1662824044243\\
11	22.1662824021965\\
12	8.34261029016324\\
13	8.34261029014876\\
14	3.03560990785045\\
15	3.03560990785033\\
16	1.34526383823868\\
17	1.34526383823868\\
18	0.484013605965423\\
19	0.484013605965423\\
20	0.276159825870787\\
21	0.276159825870787\\
22	0.128185770929315\\
23	0.128185770929315\\
24	0.0462067057407778\\
25	0.0462067057407777\\
26	0.0289707590643731\\
27	0.0289707590643731\\
28	0.0130058547326234\\
29	0.0130058547326234\\
30	0.00368459450111097\\
};
\addplot [color=mycolor2, forget plot]
  table[row sep=crcr]{%
1	0\\
2	412.935543317943\\
3	329.215137349521\\
4	53.6510157861809\\
5	53.6494408291391\\
6	15.8323674480208\\
7	15.8323619471701\\
8	4.82316263306888\\
9	4.82316239231265\\
10	2.21564594512405\\
11	2.21564586231663\\
12	0.833745445393284\\
13	0.83374542143523\\
14	0.303328559922513\\
15	0.303328549930488\\
16	0.134403192150567\\
17	0.134403188071372\\
18	0.048349108683576\\
19	0.0483491069798889\\
20	0.0275824386245548\\
21	0.0275824371338839\\
22	0.0127996381273596\\
23	0.0127996377670492\\
24	0.00461298890548032\\
25	0.00461298871873984\\
26	0.00289185583077754\\
27	0.00289185566272783\\
28	0.00129784724564763\\
29	0.00129784721240152\\
30	0.000367614458711075\\
};
\addplot [color=mycolor3, forget plot]
  table[row sep=crcr]{%
1	0\\
2	41.2935543317943\\
3	12.5592596271902\\
4	3.20274237689931\\
5	0.873643176312937\\
6	0.40299791375279\\
7	0.123855362914239\\
8	0.0570279211072458\\
9	0.0145504630244986\\
10	0.0104551183928547\\
11	0.00277652262339494\\
12	0.00161672243658461\\
13	0.000463776942426701\\
14	0.000256123539453024\\
15	6.51966694718173e-05\\
16	4.51260569620261e-05\\
};
\addplot [color=mycolor4, forget plot]
  table[row sep=crcr]{%
1	0\\
2	4.12935543317943\\
3	4.06047414977271e-08\\
4	7.87523717300535e-05\\
5	3.01540844615966e-13\\
};
\addplot [color=mycolor5, forget plot]
  table[row sep=crcr]{%
1	0\\
2	0.412935543317943\\
3	4.13671081981541e-16\\
};
\addplot [color=mycolor1, forget plot]
  table[row sep=crcr]{%
1	0\\
2	3831.80371498217\\
3	3055.76207199171\\
4	482.474032063291\\
5	482.462258992813\\
6	150.629211112038\\
7	150.6291777779\\
8	46.5385638993899\\
9	46.5385637980606\\
10	18.9308453344281\\
11	18.9308453330618\\
12	7.38241364605507\\
13	7.38241364604499\\
14	2.13874996538933\\
15	2.1387499653893\\
16	0.799854093672286\\
17	0.799854093672285\\
18	0.302401067886554\\
19	0.302401067886554\\
20	0.146988052114889\\
21	0.146988052114889\\
22	0.0714413668815339\\
23	0.0714413668815339\\
24	0.025398561546331\\
25	0.025398561546331\\
26	0.0156242317246534\\
27	0.0156242317246534\\
28	0.00815123368258414\\
29	0.00815123368258414\\
30	0.00310819177832436\\
};
\addplot [color=mycolor2, forget plot]
  table[row sep=crcr]{%
1	0\\
2	383.180371498218\\
3	305.55145426018\\
4	48.2441663376207\\
5	48.2429099889447\\
6	15.0600883282853\\
7	15.0600823305315\\
8	4.65238367684208\\
9	4.65238346447935\\
10	1.89224388762915\\
11	1.89224381793794\\
12	0.737792467443459\\
13	0.737792447335849\\
14	0.213712882808996\\
15	0.213712877486046\\
16	0.0799152875041464\\
17	0.0799152849597496\\
18	0.0302089750519511\\
19	0.0302089740961347\\
20	0.0146816132809234\\
21	0.0146816126232702\\
22	0.00713432401639561\\
23	0.00713432380377519\\
24	0.00253589765232605\\
25	0.00253589755439057\\
26	0.00155977484493652\\
27	0.00155977474386083\\
28	0.000813500628081022\\
29	0.00081350060331824\\
30	0.000310130712041444\\
};
\addplot [color=mycolor3, forget plot]
  table[row sep=crcr]{%
1	0\\
2	38.3180371498218\\
3	11.6622554008321\\
4	2.88086959208619\\
5	0.783558651143805\\
6	0.382138489507497\\
7	0.113032984419169\\
8	0.0532515542957148\\
9	0.0148978743877275\\
10	0.0092963162163118\\
11	0.00237810120872269\\
12	0.00145226246706882\\
13	0.000435108549454739\\
14	0.000190066806665846\\
15	5.73012875957189e-05\\
16	3.23063159259633e-05\\
};
\addplot [color=mycolor4, forget plot]
  table[row sep=crcr]{%
1	0\\
2	3.83180371498217\\
3	3.77166416324487e-08\\
4	7.09077879816516e-05\\
5	2.70079980830157e-13\\
};
\addplot [color=mycolor5, forget plot]
  table[row sep=crcr]{%
1	0\\
2	0.383180371498218\\
3	5.34033172918787e-16\\
};
\end{axis}
\end{tikzpicture}%
	\end{subfigure}
	\caption{Convergence of the norms $\norm{r_1}_{P_V^{-1}}$ (left) and $\norm{r_2}_{P_Q^{-1}}$ (right) of the residuals $r_1 \in V^*$ and $r_2 \in Q^*$ for the three discretizations of the elasticity problem (\cref{subsec:experiments_elasticity}) with preconditioner \eqref{eq:elasticity_preconditioner_blocks}.}
	\label{fig:residuals_Elasticity3Dpressure}
\end{figure}
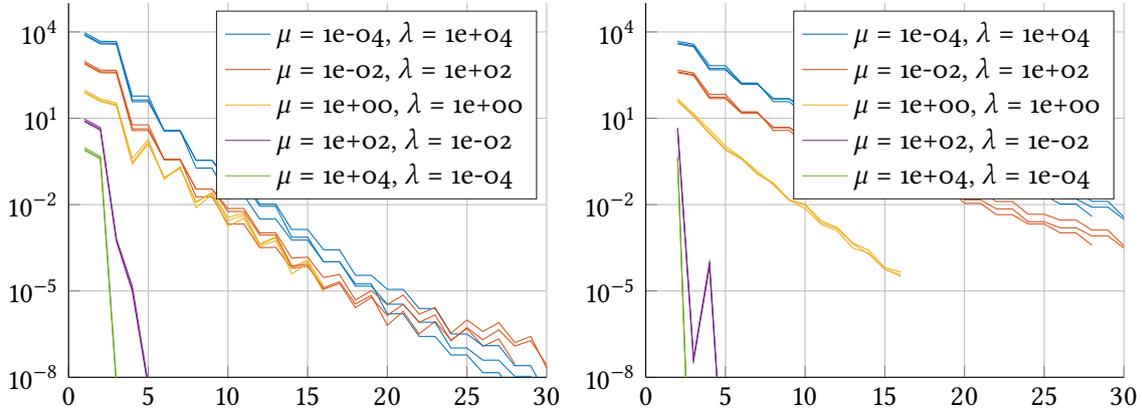

This problem is of the form \eqref{eq:elasticity_forms} and we employ the dimensionally consistent preconditioner \eqref{eq:elasticity_preconditioner_blocks}, which is robust w.r.t.\ the parameters $\mu$ and $\lambda$ as well as the mesh size.
The convergence results reported in \cref{fig:residuals_Elasticity3Dpressure,tab:convergence_Elasticity3Dpressure} confirm this.
In order to limit the amount of information, the convergence plot shows only five out of the 25~different parameter combinations.

\makeatletter
\begin{table}[htpb]
	\centering
	\begin{tabular}{@{}rrrrrrrrrr@{}}
  \toprule
           level &       $\dim(V)$ &       $\dim(Q)$ &                         $\lambda$ in \si{\newton\per\meter\squared} & $\mu = 10^{-4}$ & $\mu = 10^{-2}$ &  $\mu = 10^{0}$ &  $\mu = 10^{2}$ &  $\mu = 10^{4}$ & \si{\newton\per\meter\squared} \\
  \midrule
  \texttt{    1} & \texttt{  3075} & \texttt{   189} & \texttt{                                                     1e-04} & \texttt{    16} & \texttt{     7} & \texttt{     5} & \texttt{     3} & \texttt{     3} &                                \\
                 &                 &                 & \texttt{                                                     1e-02} & \texttt{    26} & \texttt{    16} & \texttt{     7} & \texttt{     5} & \texttt{     3} &                                \\
                 &                 &                 & \texttt{                                                     1e+00} & \texttt{    28} & \texttt{    26} & \texttt{    16} & \texttt{     7} & \texttt{     5} &                                \\
                 &                 &                 & \texttt{                                                     1e+02} & \texttt{    28} & \texttt{    28} & \texttt{    26} & \texttt{    16} & \texttt{     7} &                                \\
                 &                 &                 & \texttt{                                                     1e+04} & \texttt{    28} & \texttt{    28} & \texttt{    28} & \texttt{    26} & \texttt{    16} &                                \\
  \midrule
  \texttt{    2} & \texttt{ 19683} & \texttt{  1025} & \texttt{                                                     1e-04} & \texttt{    16} & \texttt{     7} & \texttt{     5} & \texttt{     3} & \texttt{     3} &                                \\
                 &                 &                 & \texttt{                                                     1e-02} & \texttt{    30} & \texttt{    16} & \texttt{     7} & \texttt{     5} & \texttt{     3} &                                \\
                 &                 &                 & \texttt{                                                     1e+00} & \texttt{    30} & \texttt{    30} & \texttt{    16} & \texttt{     7} & \texttt{     5} &                                \\
                 &                 &                 & \texttt{                                                     1e+02} & \texttt{    30} & \texttt{    30} & \texttt{    30} & \texttt{    16} & \texttt{     7} &                                \\
                 &                 &                 & \texttt{                                                     1e+04} & \texttt{    30} & \texttt{    30} & \texttt{    30} & \texttt{    30} & \texttt{    16} &                                \\
  \midrule
  \texttt{    3} & \texttt{139587} & \texttt{  6561} & \texttt{                                                     1e-04} & \texttt{    16} & \texttt{     7} & \texttt{     5} & \texttt{     3} & \texttt{     3} &                                \\
                 &                 &                 & \texttt{                                                     1e-02} & \texttt{    28} & \texttt{    16} & \texttt{     7} & \texttt{     5} & \texttt{     3} &                                \\
                 &                 &                 & \texttt{                                                     1e+00} & \texttt{    30} & \texttt{    28} & \texttt{    16} & \texttt{     7} & \texttt{     5} &                                \\
                 &                 &                 & \texttt{                                                     1e+02} & \texttt{    30} & \texttt{    30} & \texttt{    28} & \texttt{    16} & \texttt{     7} &                                \\
                 &                 &                 & \texttt{                                                     1e+04} & \texttt{    30} & \texttt{    30} & \texttt{    30} & \texttt{    28} & \texttt{    16} &                                \\
  \bottomrule
\end{tabular}

	\caption{Iteration numbers for the elasticity problem (\cref{subsec:experiments_elasticity}) required to reach a relative reduction by $10^{-6}$ of the initial residual norm \eqref{eq:total_residual_norm}.}
	\label{tab:convergence_Elasticity3Dpressure}
\end{table}
\makeatother

\subsection{Distributed Optimal Control of the Poisson Equation}
\label{subsec:experiments_distributed_optimal_control_Poisson}

In this section we consider the optimal control problem from \cref{subsec:distributed_optimal_control_Poisson} on $\Omega = (0,1)^3~\si{\meter\cubed}$ and with $u_d(x) = x_1~\si{\kelvin\per\meter}$ and homogeneous Dirichlet boundary conditions on all of $\Gamma$.
We discretized the problem using continuous, piecewise linear finite elements for both the state~$u$ and adjoint state~$p$.
This problem is of the form \eqref{eq:optimal_control_problem_Poisson_forms} and we employ the dimensionally consistent preconditioner \eqref{eq:optimal_control_problem_Poisson_preconditioner_particular}, which is robust w.r.t.\ \emph{all} problem parameters $\alpha$, $\beta$ and $\kappa$, as well as the mesh size.
The convergence results reported in \cref{fig:residuals_OptimalControl3DPoisson,tab:convergence_OptimalControl3DPoisson} confirm this.
While our experiments comprise the cases $\beta \in \{10^{-4}, 10^0, 10^4\}~\si{\obj\per\kelvin\squared\per\meter\cubed}$ we report only results for $\beta = 1~\si{\obj\per\kelvin\squared\per\meter\cubed}$ and $\beta = 10^4~\si{\obj\per\kelvin\squared\per\meter\cubed}$.
The maximum number of iterations was 75 and it was obtained when $\alpha = 10^{-4}~\si{\obj\per\watt\squared\meter\cubed}$, $\beta = 10^{-4}~\si{\obj\per\kelvin\squared\per\meter\cubed}$, $\kappa = 10^{-4}~\si{\watt\per\meter\per\kelvin}$ and on the finest mesh level.

Interestingly, the second residual is identically equal to zero for each instance of the problem and throughout the entire \minres iteration. 
Therefore, \cref{fig:residuals_OptimalControl3DPoisson} shows only the evolution of the norm of the first residual $\norm{r_1}_{P_V^{-1}}$.

\begin{figure}[htp]
	\begin{subfigure}[b]{0.42\textwidth}
%
\definecolor{mycolor1}{rgb}{0.00000,0.44700,0.74100}%
\definecolor{mycolor2}{rgb}{0.85000,0.32500,0.09800}%
\definecolor{mycolor3}{rgb}{0.92900,0.69400,0.12500}%
\begin{tikzpicture}

\begin{axis}[%
width=0.951\linewidth,
height=0.75\linewidth,
at={(0\linewidth,0\linewidth)},
scale only axis,
xmin=0,
xmax=35,
ymode=log,
ymin=1e-08,
ymax=1,
yminorticks=true,
axis background/.style={fill=white},
axis x line*=bottom,
axis y line*=left,
xmajorgrids,
ymajorgrids,
yminorgrids,
legend style={legend cell align=left, align=left, draw=white!15!black}
]
\addplot [color=mycolor1]
  table[row sep=crcr]{%
1	0.00623937716827497\\
2	3.48617570092551e-07\\
3	1.08031476400581e-11\\
};
\addlegendentry{$\alpha\text{ = 1e-04, }\kappa\text{ = 1e+04}$}

\addplot [color=mycolor2]
  table[row sep=crcr]{%
1	0.0616487462235321\\
2	2.64715224064502e-17\\
};
\addlegendentry{$\alpha\text{ = 1e+00, }\kappa\text{ = 1e+00}$}

\addplot [color=mycolor3]
  table[row sep=crcr]{%
1	0.343400089203619\\
2	0.0563382974489532\\
3	0.0127049105719221\\
4	0.00271940926620801\\
5	0.000783882094857128\\
6	0.000250089591724577\\
7	7.60730620180193e-05\\
8	1.97493780089764e-05\\
9	4.19135477204447e-06\\
10	5.84821966710425e-07\\
11	3.80292920844619e-08\\
};
\addlegendentry{$\alpha\text{ = 1e+04, }\kappa\text{ = 1e-04}$}

\addplot [color=mycolor1, forget plot]
  table[row sep=crcr]{%
1	0.007287694302935\\
2	6.26557601193499e-07\\
3	3.11658204124085e-11\\
};
\addplot [color=mycolor2, forget plot]
  table[row sep=crcr]{%
1	0.0719009540726281\\
2	7.99579129316044e-17\\
};
\addplot [color=mycolor3, forget plot]
  table[row sep=crcr]{%
1	0.394272302812209\\
2	0.0952536396540689\\
3	0.0353779450380038\\
4	0.0149769498650133\\
5	0.00667440354888286\\
6	0.0034748507336312\\
7	0.00183147867598805\\
8	0.000980564416138201\\
9	0.000499527071673112\\
10	0.000260972284659741\\
11	0.000138303408139624\\
12	7.22119574369394e-05\\
13	3.20831062419884e-05\\
14	1.54809792433627e-05\\
15	8.270361598014e-06\\
16	3.25972869572226e-06\\
17	1.34619411229245e-06\\
18	6.71649967688575e-07\\
19	2.405910671526e-07\\
};
\addplot [color=mycolor1, forget plot]
  table[row sep=crcr]{%
1	0.0075006457665438\\
2	6.88010032447003e-07\\
3	3.75948910251376e-11\\
};
\addplot [color=mycolor2, forget plot]
  table[row sep=crcr]{%
1	0.0739687165632564\\
2	2.014559882563e-16\\
};
\addplot [color=mycolor3, forget plot]
  table[row sep=crcr]{%
1	0.402657976479963\\
2	0.102314406957314\\
3	0.0458922943641205\\
4	0.0250088904992568\\
5	0.0148468559684851\\
6	0.00789905755299166\\
7	0.00445552593802691\\
8	0.00242809239081229\\
9	0.00158089289803939\\
10	0.000965997316072631\\
11	0.000621263486657872\\
12	0.000425034553586044\\
13	0.000263891457038076\\
14	0.000163640215093866\\
15	0.00010191189350136\\
16	6.38367657427953e-05\\
17	4.33683657028689e-05\\
18	2.51717513545033e-05\\
19	1.45071772757932e-05\\
20	8.56950887724379e-06\\
21	5.5748910017153e-06\\
22	3.59294499907519e-06\\
23	2.78226129487534e-06\\
24	1.83347417266285e-06\\
25	1.11573677719722e-06\\
26	6.90303097711832e-07\\
27	4.44462441626874e-07\\
28	3.68982140214134e-07\\
};
\addplot [color=mycolor1, forget plot]
  table[row sep=crcr]{%
1	0.00756232833703031\\
2	7.07073981694435e-07\\
3	3.95433404475826e-11\\
};
\addplot [color=mycolor2, forget plot]
  table[row sep=crcr]{%
1	0.0745694725914111\\
2	6.48618867848359e-16\\
};
\addplot [color=mycolor3, forget plot]
  table[row sep=crcr]{%
1	0.405502141343834\\
2	0.104849417843738\\
3	0.0489557431177151\\
4	0.0283841925532849\\
5	0.0185073720488212\\
6	0.0115511333038875\\
7	0.00754182587499367\\
8	0.00485486521185515\\
9	0.00329464696923761\\
10	0.00214046795470212\\
11	0.00133986323223071\\
12	0.000831336159593702\\
13	0.00058852240180538\\
14	0.000358615075531934\\
15	0.000245773168903\\
16	0.000166180917725821\\
17	0.000114579115573151\\
18	7.69275172622464e-05\\
19	5.27570603205875e-05\\
20	3.44611731069296e-05\\
21	2.45415738602832e-05\\
22	1.62113314987631e-05\\
23	1.0663108764366e-05\\
24	9.16924680490695e-06\\
25	5.92710842929345e-06\\
26	3.72037819207656e-06\\
27	2.23410792031579e-06\\
28	1.67453149658471e-06\\
29	1.39522032978569e-06\\
30	9.73972344697472e-07\\
31	6.38505116806433e-07\\
32	4.05630641996909e-07\\
33	3.05580910123823e-07\\
};
\end{axis}
\end{tikzpicture}%
	\end{subfigure}
	\hspace{0.05\textwidth}
	\begin{subfigure}[b]{0.42\textwidth}
%
\definecolor{mycolor1}{rgb}{0.00000,0.44700,0.74100}%
\definecolor{mycolor2}{rgb}{0.85000,0.32500,0.09800}%
\definecolor{mycolor3}{rgb}{0.92900,0.69400,0.12500}%
\begin{tikzpicture}

\begin{axis}[%
width=0.951\linewidth,
height=0.75\linewidth,
at={(0\linewidth,0\linewidth)},
scale only axis,
xmin=0,
xmax=30,
ymode=log,
ymin=1e-08,
ymax=1,
yminorticks=true,
axis background/.style={fill=white},
axis x line*=bottom,
axis y line*=left,
xmajorgrids,
ymajorgrids,
yminorgrids,
legend style={legend cell align=left, align=left, draw=white!15!black}
]
\addplot [color=mycolor1]
  table[row sep=crcr]{%
1	0.000616487462235321\\
2	2.76706458075215e-19\\
};
\addlegendentry{$\alpha\text{ = 1e-04, }\kappa\text{ = 1e+04}$}

\addplot [color=mycolor2]
  table[row sep=crcr]{%
1	0.00343400089203619\\
2	0.000563382974489532\\
3	0.000127049105719221\\
4	2.71940926620801e-05\\
5	7.83882094857128e-06\\
6	2.50089591724577e-06\\
7	7.60730620180195e-07\\
8	1.97493780089765e-07\\
9	4.19135477204448e-08\\
10	5.84821966710426e-09\\
11	3.80292920844618e-10\\
};
\addlegendentry{$\alpha\text{ = 1e+00, }\kappa\text{ = 1e+00}$}

\addplot [color=mycolor3]
  table[row sep=crcr]{%
1	0.00426952299314638\\
2	2.33830064162583e-05\\
3	2.05965791760076e-07\\
4	3.28153287989319e-09\\
};
\addlegendentry{$\alpha\text{ = 1e+04, }\kappa\text{ = 1e-04}$}

\addplot [color=mycolor1, forget plot]
  table[row sep=crcr]{%
1	0.000719009540726282\\
2	6.95930379823915e-19\\
};
\addplot [color=mycolor2, forget plot]
  table[row sep=crcr]{%
1	0.00394272302812209\\
2	0.000952536396540689\\
3	0.000353779450380038\\
4	0.000149769498650133\\
5	6.67440354888287e-05\\
6	3.4748507336312e-05\\
7	1.83147867598805e-05\\
8	9.80564416138202e-06\\
9	4.99527071673113e-06\\
10	2.60972284659742e-06\\
11	1.38303408139624e-06\\
12	7.22119574369395e-07\\
13	3.20831062419885e-07\\
14	1.54809792433628e-07\\
15	8.27036159801263e-08\\
16	3.25972869564269e-08\\
17	1.34619410271701e-08\\
18	6.71649528221436e-09\\
19	2.40509480991373e-09\\
};
\addplot [color=mycolor3, forget plot]
  table[row sep=crcr]{%
1	0.00515796808566517\\
2	8.66152608036182e-05\\
3	6.91698857486467e-06\\
4	5.80504841929037e-07\\
5	5.24124523329293e-08\\
6	4.65701858594833e-09\\
};
\addplot [color=mycolor1, forget plot]
  table[row sep=crcr]{%
1	0.000739687165632568\\
2	2.46754095943406e-18\\
};
\addplot [color=mycolor2, forget plot]
  table[row sep=crcr]{%
1	0.00402657976479963\\
2	0.00102314406957314\\
3	0.000458922943641205\\
4	0.000250088904992568\\
5	0.000148468559684851\\
6	7.89905755299166e-05\\
7	4.45552593802692e-05\\
8	2.4280923908123e-05\\
9	1.58089289803939e-05\\
10	9.65997316072629e-06\\
11	6.2126348665787e-06\\
12	4.25034553586043e-06\\
13	2.63891457038075e-06\\
14	1.63640215093865e-06\\
15	1.0191189350136e-06\\
16	6.38367657427945e-07\\
17	4.33683657028508e-07\\
18	2.51717513525893e-07\\
19	1.45071771440058e-07\\
20	8.56949786181327e-08\\
21	5.57441207341371e-08\\
22	3.56666599788472e-08\\
23	2.3509272006469e-08\\
24	1.82569971373251e-08\\
25	1.11568576634505e-08\\
26	6.89846784040706e-09\\
27	4.37173939538259e-09\\
28	3.55251571631757e-09\\
};
\addplot [color=mycolor3, forget plot]
  table[row sep=crcr]{%
1	0.00539940049176944\\
2	0.00019285439104494\\
3	2.92080760754132e-05\\
4	7.33538876113403e-06\\
5	1.85028555100912e-06\\
6	4.5295494470659e-07\\
7	1.14493909505591e-07\\
8	2.65794547542948e-08\\
9	6.79363556146344e-09\\
10	1.66880109577794e-09\\
};
\addplot [color=mycolor1, forget plot]
  table[row sep=crcr]{%
1	0.000745694725914108\\
2	5.50008400497297e-18\\
};
\addplot [color=mycolor2, forget plot]
  table[row sep=crcr]{%
1	0.00405502141343835\\
2	0.00104849417843739\\
3	0.000489557431177152\\
4	0.00028384192553285\\
5	0.000185073720488213\\
6	0.000115511333038876\\
7	7.5418258749937e-05\\
8	4.85486521185517e-05\\
9	3.29464696923763e-05\\
10	2.14046795470213e-05\\
11	1.33986323223072e-05\\
12	8.31336159593706e-06\\
13	5.88522401805382e-06\\
14	3.58615075531935e-06\\
15	2.45773168903001e-06\\
16	1.66180917725823e-06\\
17	1.14579115573198e-06\\
18	7.69275172639662e-07\\
19	5.27570604430817e-07\\
20	3.4461179582553e-07\\
21	2.45418611278891e-07\\
22	1.62256855256047e-07\\
23	1.13985099227872e-07\\
24	9.43956386810963e-08\\
25	5.92816742164288e-08\\
26	3.72036900034385e-08\\
27	2.2332464595679e-08\\
28	1.66356193391436e-08\\
29	1.3879888705961e-08\\
30	9.7385702513955e-09\\
31	6.38543265275905e-09\\
32	4.07056707008297e-09\\
33	3.16456709626458e-09\\
};
\addplot [color=mycolor3, forget plot]
  table[row sep=crcr]{%
1	0.00551487501435056\\
2	0.000341588430147331\\
3	7.43817038325767e-05\\
4	2.61416817074271e-05\\
5	1.25358873655416e-05\\
6	5.58965952031205e-06\\
7	2.76802304883845e-06\\
8	1.35185572533003e-06\\
9	6.32337681942687e-07\\
10	3.14119333043691e-07\\
11	1.45073157937913e-07\\
12	6.95269393768402e-08\\
13	3.23567509164069e-08\\
14	1.53471891335129e-08\\
15	7.60300883243407e-09\\
16	3.58187301737072e-09\\
};
\end{axis}
\end{tikzpicture}%
	\end{subfigure}
	\caption{Convergence of the norm $\norm{r_1}_{P_V^{-1}}$ of the residual $r_1 \in V^*$ for the four discretizations of the optimal control problem for the Poisson equation (\cref{subsec:experiments_distributed_optimal_control_Poisson}) with preconditioner \eqref{eq:optimal_control_problem_Poisson_preconditioner_particular}. The second residual $r_2 \in Q^*$ is identically equal to zero and not shown. The parameter $\beta$ is fixed to $\beta = 1~\si{\obj\per\kelvin\squared\per\meter\cubed}$ (left) and $\beta = 10^4~\si{\obj\per\kelvin\squared\per\meter\cubed}$ (right).}
	\label{fig:residuals_OptimalControl3DPoisson}
\end{figure}

\begin{table}[htpb]
	\begin{subtable}[b]{\textwidth}
		\centering
		\begin{tabular}{@{}rrrrrrrr@{}}
  \toprule
           level &       $\dim(V)$ &       $\dim(Q)$ & $\kappa$ in \si{\watt\per\meter\per\kelvin} & $\alpha = 10^{-4}$ &  $\alpha = 10^{0}$ &  $\alpha = 10^{4}$ & \si{\obj\per\watt\squared\meter\cubed} \\
  \midrule
  \texttt{    1} & \texttt{   125} & \texttt{   125} & \texttt{                             1e-04} & \texttt{        4} & \texttt{        2} & \texttt{       11} &                                        \\
                 &                 &                 & \texttt{                             1e+00} & \texttt{        9} & \texttt{        2} & \texttt{        4} &                                        \\
                 &                 &                 & \texttt{                             1e+04} & \texttt{        3} & \texttt{        2} & \texttt{        2} &                                        \\
  \midrule
  \texttt{    2} & \texttt{   729} & \texttt{   729} & \texttt{                             1e-04} & \texttt{        7} & \texttt{        2} & \texttt{       19} &                                        \\
                 &                 &                 & \texttt{                             1e+00} & \texttt{       10} & \texttt{        2} & \texttt{        4} &                                        \\
                 &                 &                 & \texttt{                             1e+04} & \texttt{        3} & \texttt{        2} & \texttt{        2} &                                        \\
  \midrule
  \texttt{    3} & \texttt{  4913} & \texttt{  4913} & \texttt{                             1e-04} & \texttt{       10} & \texttt{        2} & \texttt{       28} &                                        \\
                 &                 &                 & \texttt{                             1e+00} & \texttt{       11} & \texttt{        2} & \texttt{        4} &                                        \\
                 &                 &                 & \texttt{                             1e+04} & \texttt{        3} & \texttt{        2} & \texttt{        2} &                                        \\
  \midrule
  \texttt{    4} & \texttt{ 35937} & \texttt{ 35937} & \texttt{                             1e-04} & \texttt{       18} & \texttt{        2} & \texttt{       33} &                                        \\
                 &                 &                 & \texttt{                             1e+00} & \texttt{       11} & \texttt{        2} & \texttt{        4} &                                        \\
                 &                 &                 & \texttt{                             1e+04} & \texttt{        3} & \texttt{        2} & \texttt{        2} &                                        \\
  \bottomrule
\end{tabular}

	\end{subtable}
	\begin{subtable}[b]{\textwidth}
		\centering
		\begin{tabular}{@{}rrrrrrrr@{}}
  \toprule
           level &       $\dim(V)$ &       $\dim(Q)$ & $\kappa$ in \si{\watt\per\meter\per\kelvin} & $\alpha = 10^{-4}$ &  $\alpha = 10^{0}$ &  $\alpha = 10^{4}$ & \si{\obj\per\watt\squared\meter\cubed} \\
  \midrule
  \texttt{    1} & \texttt{   125} & \texttt{   125} & \texttt{                             1e-04} & \texttt{       15} & \texttt{        9} & \texttt{        2} &                                        \\
                 &                 &                 & \texttt{                             1e+00} & \texttt{        4} & \texttt{        3} & \texttt{        2} &                                        \\
                 &                 &                 & \texttt{                             1e+04} & \texttt{        2} & \texttt{        2} & \texttt{        2} &                                        \\
  \midrule
  \texttt{    2} & \texttt{   729} & \texttt{   729} & \texttt{                             1e-04} & \texttt{       34} & \texttt{       10} & \texttt{        2} &                                        \\
                 &                 &                 & \texttt{                             1e+00} & \texttt{        4} & \texttt{        3} & \texttt{        2} &                                        \\
                 &                 &                 & \texttt{                             1e+04} & \texttt{        2} & \texttt{        2} & \texttt{        2} &                                        \\
  \midrule
  \texttt{    3} & \texttt{  4913} & \texttt{  4913} & \texttt{                             1e-04} & \texttt{       61} & \texttt{       11} & \texttt{        2} &                                        \\
                 &                 &                 & \texttt{                             1e+00} & \texttt{        4} & \texttt{        3} & \texttt{        2} &                                        \\
                 &                 &                 & \texttt{                             1e+04} & \texttt{        2} & \texttt{        2} & \texttt{        2} &                                        \\
  \midrule
  \texttt{    4} & \texttt{ 35937} & \texttt{ 35937} & \texttt{                             1e-04} & \texttt{       75} & \texttt{       11} & \texttt{        2} &                                        \\
                 &                 &                 & \texttt{                             1e+00} & \texttt{        4} & \texttt{        3} & \texttt{        2} &                                        \\
                 &                 &                 & \texttt{                             1e+04} & \texttt{        2} & \texttt{        2} & \texttt{        2} &                                        \\
  \bottomrule
\end{tabular}

	\end{subtable}
	\caption{Iteration numbers for the optimal control problem for the Poisson equation (\cref{subsec:experiments_distributed_optimal_control_Poisson}) required to reach a relative reduction by $10^{-6}$ of the initial residual norm \eqref{eq:total_residual_norm}. The parameter $\beta$ is fixed to $1$~\si{\obj\per\kelvin\squared\per\meter\cubed} (top) and $10^{-4}$~\si{\obj\per\kelvin\squared\per\meter\cubed} (bottom).}
	\label{tab:convergence_OptimalControl3DPoisson}
\end{table}

\subsection{Distributed Optimal Control of the Stokes Equation}
\label{subsec:experiments_distributed_optimal_control_Stokes}

\begin{figure}[htp]
	\begin{subfigure}[b]{0.42\textwidth}
%
\definecolor{mycolor1}{rgb}{0.00000,0.44700,0.74100}%
\definecolor{mycolor2}{rgb}{0.85000,0.32500,0.09800}%
\definecolor{mycolor3}{rgb}{0.92900,0.69400,0.12500}%
\begin{tikzpicture}

\begin{axis}[%
width=0.951\linewidth,
height=0.75\linewidth,
at={(0\linewidth,0\linewidth)},
scale only axis,
xmin=0,
xmax=35,
ymode=log,
ymin=1e-08,
ymax=1,
yminorticks=true,
axis background/.style={fill=white},
axis x line*=bottom,
axis y line*=left,
xmajorgrids,
ymajorgrids,
yminorgrids,
legend style={legend cell align=left, align=left, draw=white!15!black}
]
\addplot [color=mycolor1]
  table[row sep=crcr]{%
1	0.00751089578000306\\
2	0.00751089546383655\\
3	0.00150217908040239\\
4	0.00150217896654808\\
5	2.25309952144914e-06\\
6	2.25171608578015e-06\\
7	2.50499047862304e-07\\
8	9.83292738187849e-08\\
9	4.52837521180426e-08\\
10	7.2256687367466e-10\\
};
\addlegendentry{$\alpha\text{ = 1e-04, }\mu\text{ = 1e+04}$}

\addplot [color=mycolor2]
  table[row sep=crcr]{%
1	0.0740559283316677\\
2	0.074025987118435\\
3	0.0148124710280015\\
4	0.0148014085994404\\
5	0.00200585876421373\\
6	0.0020028957455829\\
7	0.000220171170771886\\
8	8.87600119613829e-05\\
9	4.74180546029998e-05\\
10	4.39545690192247e-05\\
11	2.33177914720165e-05\\
12	1.98172038160725e-05\\
13	1.13628169200498e-06\\
14	3.16041125516302e-07\\
15	2.3403429800231e-07\\
16	2.12337376451887e-07\\
17	8.24042973154871e-08\\
18	6.93667367046586e-08\\
19	6.30595834543749e-09\\
};
\addlegendentry{$\alpha\text{ = 1e+00, }\mu\text{ = 1e+00}$}

\addplot [color=mycolor3]
  table[row sep=crcr]{%
1	0.399675760054392\\
2	0.290385203097822\\
3	0.0486806313741074\\
4	0.033913048389886\\
5	0.0196751777222307\\
6	0.0205635216779355\\
7	0.00712140513752078\\
8	0.00996997822204462\\
9	0.00356575875914921\\
10	0.00300355363499611\\
11	0.00148375654459521\\
12	0.00148789466417722\\
13	0.000607037071258901\\
14	0.000357122728460078\\
15	0.000235622555483825\\
16	0.000245672972811564\\
17	7.84306875351709e-05\\
18	5.01455873564175e-05\\
19	2.63578363645278e-05\\
20	2.69357721864462e-05\\
21	9.5661488710686e-06\\
22	7.27796514779288e-06\\
23	4.55199748115473e-06\\
24	4.63952068375903e-06\\
25	2.15154005530307e-06\\
26	1.36061210453894e-06\\
27	6.7213260075014e-07\\
28	6.69663368044153e-07\\
29	2.52415244204741e-07\\
30	1.83586772073856e-07\\
};
\addlegendentry{$\alpha\text{ = 1e+04, }\mu\text{ = 1e-04}$}

\addplot [color=mycolor1, forget plot]
  table[row sep=crcr]{%
1	0.00757804466043454\\
2	0.00757804434345104\\
3	0.00151562920568004\\
4	0.00151562909153762\\
5	1.67115919825515e-05\\
6	1.67115743645456e-05\\
7	3.3340609294846e-06\\
8	3.33406078205154e-06\\
9	6.56993608698484e-08\\
10	6.56993597436951e-08\\
11	1.10385944291526e-09\\
};
\addplot [color=mycolor2, forget plot]
  table[row sep=crcr]{%
1	0.0747215407342616\\
2	0.0746915476153482\\
3	0.0149469251014452\\
4	0.0149358320890764\\
5	0.00212364858395714\\
6	0.00212071824034155\\
7	0.000248676709406686\\
8	0.000133011701196745\\
9	5.98959540022108e-05\\
10	6.05533333963706e-05\\
11	4.41077279827099e-05\\
12	4.40327872926947e-05\\
13	1.5385730744163e-06\\
14	5.66364496489253e-07\\
15	3.52709503510779e-07\\
16	3.32785842227443e-07\\
17	1.26024215911252e-07\\
18	1.28479504264749e-07\\
19	1.03390458770761e-08\\
};
\addplot [color=mycolor3, forget plot]
  table[row sep=crcr]{%
1	0.405962083883989\\
2	0.299208246038815\\
3	0.0542700086775693\\
4	0.0357503983399989\\
5	0.024226848854271\\
6	0.0244287696587932\\
7	0.0076232980708994\\
8	0.01067689009982\\
9	0.00489693045426498\\
10	0.0049418399265248\\
11	0.00202857816336506\\
12	0.00191247724361483\\
13	0.000988942646794895\\
14	0.000793269648336384\\
15	0.000328803406459393\\
16	0.000324434170792504\\
17	0.000147509204616331\\
18	0.00013719080968491\\
19	4.54358037612033e-05\\
20	4.42385039456078e-05\\
21	1.94761601462005e-05\\
22	1.7799629728016e-05\\
23	8.1524730606067e-06\\
24	7.79786121477633e-06\\
25	3.65257120283699e-06\\
26	3.27026123612119e-06\\
27	1.36936990142738e-06\\
28	1.34220944616384e-06\\
29	6.02027611394559e-07\\
30	4.50183401863764e-07\\
31	2.23087668364088e-07\\
32	2.03972120951583e-07\\
};
\addplot [color=mycolor1, forget plot]
  table[row sep=crcr]{%
1	0.00758330486229702\\
2	0.00758330454535241\\
3	0.0015166625216026\\
4	0.00151666240746705\\
5	5.12019217625012e-06\\
6	5.11988009420098e-06\\
7	2.19219839228861e-06\\
8	2.19210777984384e-06\\
9	6.60235876115829e-08\\
10	6.60235864705415e-08\\
11	1.10648434288796e-09\\
};
\addplot [color=mycolor2, forget plot]
  table[row sep=crcr]{%
1	0.07477385033579\\
2	0.0747438627201339\\
3	0.0149572877882246\\
4	0.0149461957764504\\
5	0.00212574736572884\\
6	0.00212279863811518\\
7	0.000232478226504171\\
8	9.76024946181694e-05\\
9	5.33575781586189e-05\\
10	5.2043817267801e-05\\
11	3.01227246447879e-05\\
12	2.78240204457383e-05\\
13	1.40554593737421e-06\\
14	5.28131909141445e-07\\
15	3.3573183975528e-07\\
16	3.08568421545451e-07\\
17	1.23144397017457e-07\\
18	1.26752851833501e-07\\
19	1.06373939376624e-08\\
};
\addplot [color=mycolor3, forget plot]
  table[row sep=crcr]{%
1	0.406486511171613\\
2	0.299972833023571\\
3	0.0548678059510683\\
4	0.0359751101035927\\
5	0.0252513284053981\\
6	0.0251283115339676\\
7	0.0081132881405856\\
8	0.0108839042732656\\
9	0.00551603122258803\\
10	0.00570811136818313\\
11	0.00220349486625525\\
12	0.00194659169915301\\
13	0.00112774470821142\\
14	0.000984617315253866\\
15	0.000377622333603274\\
16	0.000345008659119952\\
17	0.00017800935995711\\
18	0.000170885276156161\\
19	5.95112798568502e-05\\
20	5.37663211274839e-05\\
21	2.5254829633707e-05\\
22	2.45393617911652e-05\\
23	1.02816902664787e-05\\
24	9.12097892918076e-06\\
25	5.02350911266245e-06\\
26	4.82704657010706e-06\\
27	1.86073524742504e-06\\
28	1.71788767454581e-06\\
29	7.75920025104332e-07\\
30	6.67768390787266e-07\\
31	2.99344029162564e-07\\
32	2.44505185772655e-07\\
33	1.33065991085252e-07\\
};
\end{axis}
\end{tikzpicture}%
	\end{subfigure}
	\hspace{0.05\textwidth}
	\begin{subfigure}[b]{0.42\textwidth}
%
\definecolor{mycolor1}{rgb}{0.00000,0.44700,0.74100}%
\definecolor{mycolor2}{rgb}{0.85000,0.32500,0.09800}%
\definecolor{mycolor3}{rgb}{0.92900,0.69400,0.12500}%
\begin{tikzpicture}

\begin{axis}[%
width=0.951\linewidth,
height=0.75\linewidth,
at={(0\linewidth,0\linewidth)},
scale only axis,
xmin=0,
xmax=35,
ymode=log,
ymin=1e-08,
ymax=1,
yminorticks=true,
axis background/.style={fill=white},
axis x line*=bottom,
axis y line*=left,
xmajorgrids,
ymajorgrids,
yminorgrids,
legend style={legend cell align=left, align=left, draw=white!15!black}
]
\addplot [color=mycolor1]
  table[row sep=crcr]{%
1	0\\
2	0\\
3	0.00300435804377816\\
4	0.00300435800582377\\
5	1.10919085333116e-06\\
6	1.10919085491429e-06\\
7	1.78183666626688e-07\\
8	1.6455036226195e-07\\
9	6.98623690842819e-08\\
10	3.49708832663001e-10\\
};
\addlegendentry{$\alpha\text{ = 1e-04, }\mu\text{ = 1e+04}$}

\addplot [color=mycolor2]
  table[row sep=crcr]{%
1	0\\
2	0\\
3	0.0296038230003589\\
4	0.0296000973377498\\
5	0.00100278682478061\\
6	0.00100279150842143\\
7	0.000163885425381573\\
8	0.000151495235405124\\
9	6.85495043771356e-05\\
10	5.23829864296533e-05\\
11	2.24460994673961e-05\\
12	9.73431442986603e-06\\
13	3.74378524638255e-07\\
14	3.79806103621371e-07\\
15	2.78072296760761e-07\\
16	2.78820083621934e-07\\
17	1.31764801460426e-07\\
18	4.64513152666079e-08\\
19	3.53109554166566e-09\\
};
\addlegendentry{$\alpha\text{ = 1e+00, }\mu\text{ = 1e+00}$}

\addplot [color=mycolor3]
  table[row sep=crcr]{%
1	0\\
2	0\\
3	0.0975069727422016\\
4	0.0710058945813276\\
5	0.0175325063908976\\
6	0.020521204884734\\
7	0.00827713225907892\\
8	0.00473590059447967\\
9	0.00242263828265599\\
10	0.00208300690968529\\
11	0.00127302934729163\\
12	0.00126006182781057\\
13	0.000605346513038336\\
14	0.000445373351815466\\
15	0.000169278687768037\\
16	0.000157489711894913\\
17	4.87989658526005e-05\\
18	3.66191957115017e-05\\
19	2.33002722003683e-05\\
20	2.08797767683465e-05\\
21	8.38130413445566e-06\\
22	6.89710947424085e-06\\
23	3.78659450345963e-06\\
24	3.72448710057691e-06\\
25	1.65602641843639e-06\\
26	1.00142496179532e-06\\
27	6.025393333131e-07\\
28	6.08014135793763e-07\\
29	2.07254918839773e-07\\
30	1.58844207283847e-07\\
};
\addlegendentry{$\alpha\text{ = 1e+04, }\mu\text{ = 1e-04}$}

\addplot [color=mycolor1, forget plot]
  table[row sep=crcr]{%
1	0\\
2	0\\
3	0.00303121477368093\\
4	0.00303121473563022\\
5	1.17593751989144e-06\\
6	1.17593750716953e-06\\
7	2.0250792014768e-06\\
8	2.02507914295647e-06\\
9	1.05106413010754e-07\\
10	1.05106413028287e-07\\
11	6.87184531159558e-10\\
};
\addplot [color=mycolor2, forget plot]
  table[row sep=crcr]{%
1	0\\
2	0\\
3	0.0298697213598217\\
4	0.0298660000123909\\
5	0.00106091600876371\\
6	0.00106092276607434\\
7	0.000173975239881477\\
8	0.000162175040124326\\
9	8.35529992406136e-05\\
10	8.10244142628269e-05\\
11	4.21457892559001e-05\\
12	4.06005785510216e-05\\
13	4.58315754859138e-07\\
14	4.64785642407386e-07\\
15	3.41650731713014e-07\\
16	3.37601942306075e-07\\
17	1.89269073576643e-07\\
18	1.81299746408665e-07\\
19	5.55124359736747e-09\\
};
\addplot [color=mycolor3, forget plot]
  table[row sep=crcr]{%
1	0\\
2	0\\
3	0.100906554490553\\
4	0.0756917662904331\\
5	0.0205187658901344\\
6	0.0224077969385743\\
7	0.00912825971288962\\
8	0.00538258722286834\\
9	0.00327394008732835\\
10	0.00304238297573802\\
11	0.00160831126084184\\
12	0.00158040769842465\\
13	0.000816320128130498\\
14	0.000765676357223475\\
15	0.000446869459177563\\
16	0.000349287560167319\\
17	7.9374533950347e-05\\
18	6.94157036638813e-05\\
19	4.86741265425822e-05\\
20	3.9773529853363e-05\\
21	1.81877971516685e-05\\
22	1.7858198462336e-05\\
23	8.21165067067381e-06\\
24	7.48778557881446e-06\\
25	2.6614747191217e-06\\
26	2.41300787931257e-06\\
27	1.43313357565542e-06\\
28	1.05261473611874e-06\\
29	4.79310075315036e-07\\
30	4.38914734081689e-07\\
31	2.28898820984668e-07\\
32	2.10185250268981e-07\\
};
\addplot [color=mycolor1, forget plot]
  table[row sep=crcr]{%
1	0\\
2	0\\
3	0.00303332152676928\\
4	0.00303332148872269\\
5	1.17390621466288e-06\\
6	1.1739062153491e-06\\
7	8.43119056437795e-07\\
8	8.43059199893613e-07\\
9	1.0549253461162e-07\\
10	1.05492534838964e-07\\
11	6.90083146257637e-10\\
};
\addplot [color=mycolor2, forget plot]
  table[row sep=crcr]{%
1	0\\
2	0\\
3	0.0298906522351316\\
4	0.0298869323679899\\
5	0.0010645668018475\\
6	0.00106457399219325\\
7	0.000173396852534648\\
8	0.000160700175919525\\
9	7.50648258201975e-05\\
10	6.33753579760241e-05\\
11	2.8987554950497e-05\\
12	1.73081617011032e-05\\
13	4.64083214023578e-07\\
14	4.69966416088901e-07\\
15	3.419822019943e-07\\
16	3.38332818973855e-07\\
17	1.83375122037759e-07\\
18	1.68190327500303e-07\\
19	6.04679624517303e-09\\
};
\addplot [color=mycolor3, forget plot]
  table[row sep=crcr]{%
1	0\\
2	0\\
3	0.10114763417654\\
4	0.0761011371742578\\
5	0.0215880159346964\\
6	0.0232902864524041\\
7	0.00923501901831018\\
8	0.0058143895967047\\
9	0.00353915710803863\\
10	0.00328719710828813\\
11	0.00179322859864801\\
12	0.00170219247954337\\
13	0.000883226683844474\\
14	0.000866739723684641\\
15	0.000511472336132032\\
16	0.000375984038907499\\
17	0.000104977736757131\\
18	9.69730230763795e-05\\
19	5.56698525182954e-05\\
20	4.30194797106171e-05\\
21	2.37778395151164e-05\\
22	2.36722345688905e-05\\
23	1.07475059330129e-05\\
24	9.32404711539447e-06\\
25	3.91173086662477e-06\\
26	3.77094459473995e-06\\
27	1.91593814788401e-06\\
28	1.31941964267358e-06\\
29	6.49704871796066e-07\\
30	6.13814687638233e-07\\
31	3.19031257991045e-07\\
32	2.77613188051585e-07\\
33	9.95170765284824e-08\\
};
\end{axis}
\end{tikzpicture}%
	\end{subfigure}
	\\
	\begin{subfigure}[b]{0.42\textwidth}
%
\definecolor{mycolor1}{rgb}{0.00000,0.44700,0.74100}%
\definecolor{mycolor2}{rgb}{0.85000,0.32500,0.09800}%
\definecolor{mycolor3}{rgb}{0.92900,0.69400,0.12500}%
\begin{tikzpicture}

\begin{axis}[%
width=0.951\linewidth,
height=0.75\linewidth,
at={(0\linewidth,0\linewidth)},
scale only axis,
xmin=0,
xmax=35,
ymode=log,
ymin=1e-08,
ymax=1,
yminorticks=true,
axis background/.style={fill=white},
axis x line*=bottom,
axis y line*=left,
xmajorgrids,
ymajorgrids,
yminorgrids,
legend style={legend cell align=left, align=left, draw=white!15!black}
]
\addplot [color=mycolor1]
  table[row sep=crcr]{%
1	0\\
2	1.08949633997381e-06\\
3	2.18000886508963e-07\\
4	4.35965582040331e-07\\
5	3.99308019708203e-07\\
6	3.96888058229344e-07\\
7	2.43327641264096e-07\\
8	5.98396800671642e-08\\
9	2.89284106939621e-08\\
10	1.02793754061847e-09\\
};
\addlegendentry{$\alpha\text{ = 1e-04, }\mu\text{ = 1e+04}$}

\addplot [color=mycolor2]
  table[row sep=crcr]{%
1	0\\
2	0.00103765108817605\\
3	0.000217170436754989\\
4	0.00043027871649381\\
5	0.000380076385316273\\
6	0.000380140729233425\\
7	0.000234562806034163\\
8	7.05941570206916e-05\\
9	2.79290951542812e-05\\
10	3.84653731140808e-05\\
11	2.72925048166363e-05\\
12	1.98807289944638e-05\\
13	1.18062495321385e-06\\
14	4.1256998552859e-07\\
15	1.44015893800368e-07\\
16	9.46541862477516e-08\\
17	6.07521300555833e-08\\
18	6.85951195429859e-08\\
19	3.09965510620572e-09\\
};
\addlegendentry{$\alpha\text{ = 1e+00, }\mu\text{ = 1e+00}$}

\addplot [color=mycolor3]
  table[row sep=crcr]{%
1	0\\
2	0.0613716777522341\\
3	0.0496181877069992\\
4	0.0316105891060811\\
5	0.0381280726006939\\
6	0.0318805946935575\\
7	0.0145509536100996\\
8	0.00952166511723238\\
9	0.00460928358745421\\
10	0.00314986560678406\\
11	0.00201594902053693\\
12	0.00118229455689205\\
13	0.000640046370730501\\
14	0.000444155203677708\\
15	0.000196559707204885\\
16	0.000178053430615661\\
17	0.0001008454382957\\
18	5.3608572655086e-05\\
19	2.69056329141802e-05\\
20	2.18716812361902e-05\\
21	1.32022423169789e-05\\
22	8.60988329338854e-06\\
23	4.14368228215364e-06\\
24	3.71757603221238e-06\\
25	2.31477856080916e-06\\
26	1.15654764963579e-06\\
27	5.56591186064349e-07\\
28	5.27923162345639e-07\\
29	3.5872204219827e-07\\
30	1.97106598276723e-07\\
};
\addlegendentry{$\alpha\text{ = 1e+04, }\mu\text{ = 1e-04}$}

\addplot [color=mycolor1, forget plot]
  table[row sep=crcr]{%
1	0\\
2	1.09577050461776e-06\\
3	2.19253757841489e-07\\
4	4.38491338610603e-07\\
5	4.16012870244699e-07\\
6	4.16004194415088e-07\\
7	3.73966290911876e-08\\
8	3.74032310686653e-08\\
9	2.18589041637248e-11\\
10	2.35008287661631e-11\\
11	7.28748136387374e-11\\
};
\addplot [color=mycolor2, forget plot]
  table[row sep=crcr]{%
1	0\\
2	0.00104326568856488\\
3	0.000218091095517827\\
4	0.000432954128411958\\
5	0.00039565265229575\\
6	0.00039573536851145\\
7	0.000244413076346313\\
8	0.000114075921546892\\
9	2.37862861431531e-05\\
10	2.79304269856172e-05\\
11	2.1299599727317e-05\\
12	2.19750257553563e-05\\
13	1.62367425235819e-06\\
14	7.03170475752108e-07\\
15	2.21166910772403e-07\\
16	1.87452949069039e-07\\
17	5.73322640348538e-08\\
18	6.71970403290179e-08\\
19	5.82258182438928e-09\\
};
\addplot [color=mycolor3, forget plot]
  table[row sep=crcr]{%
1	0\\
2	0.0640140751731157\\
3	0.0485177406507589\\
4	0.033768934820037\\
5	0.0388282503348651\\
6	0.0328868259794086\\
7	0.0200302154114567\\
8	0.0112241964283783\\
9	0.00526747445937297\\
10	0.00434220767243863\\
11	0.00288078959719069\\
12	0.00165941201988782\\
13	0.000802600526011402\\
14	0.000665185868252193\\
15	0.00044763666770832\\
16	0.000327119438021616\\
17	0.000146554493306086\\
18	0.000111424717420087\\
19	6.14026680533627e-05\\
20	4.19419676658724e-05\\
21	1.89432757070176e-05\\
22	1.65598473640191e-05\\
23	1.10828225995943e-05\\
24	7.89371778808801e-06\\
25	3.39262738487422e-06\\
26	2.84626960124001e-06\\
27	1.79603732615556e-06\\
28	1.16908379711056e-06\\
29	5.87820304589315e-07\\
30	4.42213457045926e-07\\
31	2.7190525964571e-07\\
32	1.92608807268832e-07\\
};
\addplot [color=mycolor1, forget plot]
  table[row sep=crcr]{%
1	0\\
2	1.09608235950129e-06\\
3	2.19315933282356e-07\\
4	4.38605944984595e-07\\
5	4.17754429410664e-07\\
6	4.17209976291605e-07\\
7	1.78258499674463e-07\\
8	1.78480237534649e-07\\
9	2.70222756094874e-10\\
10	2.70296702972535e-10\\
11	7.52504219931871e-11\\
};
\addplot [color=mycolor2, forget plot]
  table[row sep=crcr]{%
1	0\\
2	0.00104354079073761\\
3	0.000218126558758761\\
4	0.000433076019172892\\
5	0.000397301623076299\\
6	0.000397379610527814\\
7	0.000245582268568432\\
8	7.9054430518034e-05\\
9	2.85796940411963e-05\\
10	3.85442035760392e-05\\
11	2.95212733156934e-05\\
12	2.55654242485488e-05\\
13	1.45504098854792e-06\\
14	6.56298269354348e-07\\
15	2.27161927657887e-07\\
16	1.77066394121733e-07\\
17	6.26256901776499e-08\\
18	7.73084646323063e-08\\
19	4.94259190945614e-09\\
};
\addplot [color=mycolor3, forget plot]
  table[row sep=crcr]{%
1	0\\
2	0.0642393386713726\\
3	0.0483958708901898\\
4	0.0339924937980296\\
5	0.0388322892131591\\
6	0.0331648806385402\\
7	0.0218137995889545\\
8	0.0118352704377086\\
9	0.0055528539730144\\
10	0.00482185189347803\\
11	0.00321329248737066\\
12	0.00176332339521324\\
13	0.000848612355218957\\
14	0.000721398442715748\\
15	0.000538327237683849\\
16	0.000374779706550896\\
17	0.000160212521069528\\
18	0.000133564543596122\\
19	8.40606928482945e-05\\
20	5.17383293368282e-05\\
21	2.14559575947949e-05\\
22	1.96313584860565e-05\\
23	1.39930293402048e-05\\
24	9.39176335698986e-06\\
25	4.11656531936463e-06\\
26	3.64552320774112e-06\\
27	2.561027274009e-06\\
28	1.59066934540082e-06\\
29	6.86059043706309e-07\\
30	5.57708582856525e-07\\
31	3.97706164299798e-07\\
32	2.50457784047648e-07\\
33	1.11338277117703e-07\\
};
\end{axis}
\end{tikzpicture}%
	\end{subfigure}
	\hspace{0.05\textwidth}
	\begin{subfigure}[b]{0.42\textwidth}
%
\definecolor{mycolor1}{rgb}{0.00000,0.44700,0.74100}%
\definecolor{mycolor2}{rgb}{0.85000,0.32500,0.09800}%
\definecolor{mycolor3}{rgb}{0.92900,0.69400,0.12500}%
\begin{tikzpicture}

\begin{axis}[%
width=0.951\linewidth,
height=0.75\linewidth,
at={(0\linewidth,0\linewidth)},
scale only axis,
xmin=0,
xmax=35,
ymode=log,
ymin=1e-08,
ymax=1,
yminorticks=true,
axis background/.style={fill=white},
axis x line*=bottom,
axis y line*=left,
xmajorgrids,
ymajorgrids,
yminorgrids,
legend style={legend cell align=left, align=left, draw=white!15!black}
]
\addplot [color=mycolor1]
  table[row sep=crcr]{%
1	0\\
2	1.0898124973146e-06\\
3	6.90405096277673e-07\\
4	2.21670009191867e-07\\
5	1.81733871018568e-10\\
6	4.33557663307013e-08\\
7	2.65831831270554e-08\\
8	9.03621235770662e-08\\
9	4.36641915362162e-08\\
10	5.82345767375108e-10\\
};
\addlegendentry{$\alpha\text{ = 1e-04, }\mu\text{ = 1e+04}$}

\addplot [color=mycolor2]
  table[row sep=crcr]{%
1	0\\
2	0.00106755048549654\\
3	0.000666997835206845\\
4	0.000199050121223716\\
5	1.72330427624361e-05\\
6	4.09462250870974e-05\\
7	2.60804209358086e-05\\
8	8.55339471358993e-05\\
9	3.22688445251786e-05\\
10	2.23706180327434e-05\\
11	1.6187853715395e-05\\
12	2.76753911577192e-07\\
13	2.54800691174384e-08\\
14	2.85986060779557e-07\\
15	1.01229353005994e-07\\
16	1.15805465806103e-07\\
17	7.15569167453665e-08\\
18	2.60728072634575e-09\\
19	8.71354565637713e-11\\
};
\addlegendentry{$\alpha\text{ = 1e+00, }\mu\text{ = 1e+00}$}

\addplot [color=mycolor3]
  table[row sep=crcr]{%
1	0\\
2	0.165759894587424\\
3	0.0596059625931801\\
4	0.0310433679794107\\
5	0.00615975803689298\\
6	0.00944100738519161\\
7	0.00420873763219784\\
8	0.00218867104293809\\
9	0.00204023232767167\\
10	0.00172773745074181\\
11	0.000882666266348837\\
12	0.000255770695095918\\
13	0.000458612199147758\\
14	0.000263027935497096\\
15	5.60926369122183e-05\\
16	4.45493569816871e-05\\
17	3.76402791932142e-05\\
18	2.85721721000216e-05\\
19	9.8317601569358e-06\\
20	6.66772482985562e-06\\
21	6.11688502582936e-06\\
22	4.36049351540811e-06\\
23	2.04481143128546e-06\\
24	1.50725769313972e-06\\
25	1.39517414685971e-06\\
26	6.17749484719351e-07\\
27	2.28072173689074e-07\\
28	1.99057048105949e-07\\
29	1.56248261309913e-07\\
30	7.56028645275566e-08\\
};
\addlegendentry{$\alpha\text{ = 1e+04, }\mu\text{ = 1e-04}$}

\addplot [color=mycolor1, forget plot]
  table[row sep=crcr]{%
1	0\\
2	1.09608486827452e-06\\
3	6.97988334531133e-07\\
4	2.33916569728711e-07\\
5	1.93531124360322e-10\\
6	3.14725670663786e-09\\
7	5.58186179323858e-10\\
8	1.27112328648479e-10\\
9	1.5272760281538e-11\\
10	6.22785646894734e-12\\
11	1.02833806446836e-13\\
};
\addplot [color=mycolor2, forget plot]
  table[row sep=crcr]{%
1	0\\
2	0.00107321054518005\\
3	0.000674412976849603\\
4	0.000210215854096376\\
5	1.81793577053361e-05\\
6	4.18348030956498e-05\\
7	2.66900415999939e-05\\
8	8.38435089820627e-05\\
9	1.59681266797165e-05\\
10	8.14522823189189e-06\\
11	6.71794877959701e-06\\
12	3.89235305847525e-07\\
13	4.28508525488707e-08\\
14	3.60223236700889e-07\\
15	1.15027633319674e-07\\
16	1.20042793707659e-07\\
17	3.40748835709898e-08\\
18	1.02676822475311e-09\\
19	2.16556944351818e-10\\
};
\addplot [color=mycolor3, forget plot]
  table[row sep=crcr]{%
1	0\\
2	0.165106350304855\\
3	0.0649151237572185\\
4	0.0326923187754447\\
5	0.0073961408219653\\
6	0.0115446957987344\\
7	0.0060117921090715\\
8	0.00210303438624242\\
9	0.00304310190226806\\
10	0.00307512675555655\\
11	0.00150262063617866\\
12	0.000268163767264187\\
13	0.00053619635109153\\
14	0.000432529704664609\\
15	0.000211996963720393\\
16	7.36379687838752e-05\\
17	4.34936084515475e-05\\
18	5.77822617104525e-05\\
19	2.25625339434178e-05\\
20	1.30621581074496e-05\\
21	9.68120613405279e-06\\
22	9.1081273527772e-06\\
23	4.83029401766971e-06\\
24	2.67150308202595e-06\\
25	1.68379357624494e-06\\
26	1.40999355208367e-06\\
27	7.48303542665397e-07\\
28	3.16393901406338e-07\\
29	2.97675323954755e-07\\
30	2.85597805951069e-07\\
31	1.20302157134974e-07\\
32	6.89225416722823e-08\\
};
\addplot [color=mycolor1, forget plot]
  table[row sep=crcr]{%
1	0\\
2	1.0963991033277e-06\\
3	6.98437962127588e-07\\
4	2.34729174825836e-07\\
5	1.91994334098009e-10\\
6	2.11579102713628e-08\\
7	9.13505626966165e-09\\
8	1.25279626996338e-09\\
9	1.36249733712429e-11\\
10	6.25569130139697e-12\\
11	1.06075428317201e-13\\
};
\addplot [color=mycolor2, forget plot]
  table[row sep=crcr]{%
1	0\\
2	0.00107348229201937\\
3	0.000674854489875174\\
4	0.000210955483766714\\
5	1.8256557618342e-05\\
6	4.2038225556528e-05\\
7	2.68181327018964e-05\\
8	8.99028865948405e-05\\
9	3.09067735668065e-05\\
10	2.12344919162965e-05\\
11	1.66654399874732e-05\\
12	2.50105725664594e-07\\
13	3.49495304143834e-08\\
14	3.61695951084132e-07\\
15	1.26896373111728e-07\\
16	1.34429951063742e-07\\
17	4.48066793600316e-08\\
18	1.41983905044919e-09\\
19	1.60387731362406e-10\\
};
\addplot [color=mycolor3, forget plot]
  table[row sep=crcr]{%
1	0\\
2	0.165022637250249\\
3	0.0653236061819401\\
4	0.0327953738451139\\
5	0.00758117735140649\\
6	0.0116083324986073\\
7	0.00639157198106031\\
8	0.00230637750320628\\
9	0.00320578078200978\\
10	0.00326046957320336\\
11	0.00171453844526879\\
12	0.000374934456014808\\
13	0.000546337317247787\\
14	0.000452740852075992\\
15	0.000269559119493663\\
16	9.55836685677097e-05\\
17	4.78110386680083e-05\\
18	5.9269128913474e-05\\
19	3.01954997665199e-05\\
20	1.82707001166085e-05\\
21	1.09680218995196e-05\\
22	1.0497292740446e-05\\
23	6.41392523854033e-06\\
24	3.56947018143146e-06\\
25	2.04324354148397e-06\\
26	1.92327748824474e-06\\
27	1.11846673155667e-06\\
28	4.46156036913743e-07\\
29	3.67905939137835e-07\\
30	3.85759209826876e-07\\
31	1.95975740926701e-07\\
32	1.01438796258392e-07\\
33	5.10484770887113e-08\\
};
\end{axis}
\end{tikzpicture}%
	\end{subfigure}
	\caption{Convergence of the norms $\norm{r_{1,1}}_{P_{V_1}^{-1}}$ (top left), $\norm{r_{1,2}}_{P_{V_2}^{-1}}$ (top right), $\norm{r_{2,1}}_{P_{Q_1}^{-1}}$ (bottom left) and $\norm{r_{2,2}}_{P_{Q_2}^{-1}}$ (bottom right) of the residuals $r_1 \in V^*$ and $r_2 \in Q^*$ for the three discretizations of the optimal control problem for the Stokes system (\cref{subsec:experiments_distributed_optimal_control_Stokes}) with preconditioner \eqref{eq:optimal_control_problem_Stokes_preconditioner_particular}. The parameter $\beta$ is fixed to $1$~\si{\obj\per\meter\tothe{5}\second\squared}.}
	\label{fig:residuals_OptimalControl3DStokes}
\end{figure}
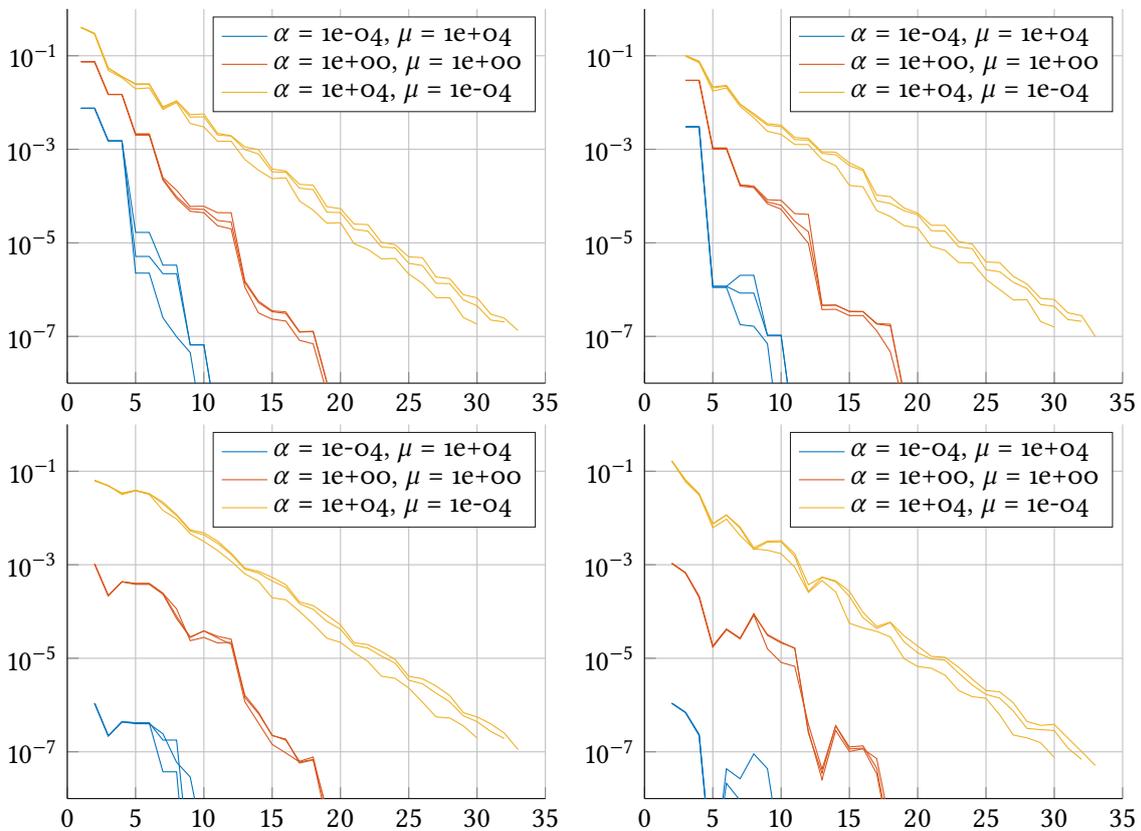

As our final example, we consider the problem from \cref{subsec:distributed_optimal_control_Stokes} on $\Omega = (0,1)^3~\si{\meter\cubed}$ and with $u_d(x) = x_1~\si{\per\second}$ and homogeneous Dirichlet boundary conditions on all of $\Gamma$.
We discretized the problem using a Taylor--Hood pair for the state velocity/pressure pair~$(u,p)$, and the same function space for the adjoint velocity/pressure pair~$(w,r)$.
This problem is of the form \eqref{eq:optimal_control_problem_Stokes_forms} and we employ the dimensionally consistent preconditioner \eqref{eq:optimal_control_problem_Stokes_preconditioner_particular}, which is robust w.r.t.\ \emph{all} problem parameters $\alpha$, $\beta$ and $\mu$, as well as the mesh size.
The convergence results reported in \cref{fig:residuals_OptimalControl3DStokes,tab:convergence_OptimalControl3DStokes} confirm this.
In order to limit the amount of information, we report only results for $\beta = 1~\si{\obj\per\meter\tothe{5}\second\squared}$ while our experiments comprise the cases $\beta \in \{10^{-4}, 10^0, 10^4\}~\si{\obj\per\meter\tothe{5}\second\squared}$.
The maximum number of iterations was 33 and it was obtained when $\alpha = 10^4~\si{\obj\per\newton\squared\meter\cubed}$, $\beta = 1~\si{\obj\per\meter\tothe{5}\second\squared}$ and $\mu = 10^{-4}~\si{\newton\second\per\meter\squared}$ and on the finest mesh level.

The application of the preconditioner \eqref{eq:optimal_control_problem_Stokes_preconditioner_particular} deserves a more detailed description in this case.
As described in the introduction to \cref{sec:numerical_experiments}, we obtain a Cholesky factorization of $\beta \, \tM + (\alpha \, \beta)^{1/2} \mu \, \tK$, which we utilize in the first block of $\tP_V$ and $\tP_Q$.
As for the second block of $\tP_V$ and $\tP_Q$, we implemented matrix-vector products with $\tD \paren[big][]{\beta \, \tM + (\alpha \, \beta)^{1/2} \mu \, \tK}^{-1} \tD^\transp$ (utilizing again the Cholesky factorization) and applied \matlab's conjugate gradient solver with default settings and no preconditioner.
Although $\tP_V$ and $\tP_Q$ then effectively become mildly nonlinear preconditioners, this did not affect the convergence of the outer \minres iteration.

\begin{table}[htpb]
	\centering
	\begin{tabular}{@{}rrrrrrrr@{}}
  \toprule
           level &             $\dim(V_1)$ &             $\dim(V_2)$ & $\mu$ in \si{\newton\second\per\meter\squared} & $\alpha = 10^{-4}$ &  $\alpha = 10^{0}$ &  $\alpha = 10^{4}$ & \si{\obj\per\newton\squared\meter\cubed} \\
  \midrule
  \texttt{    1} & \texttt{          2187} & \texttt{           125} & \texttt{                                1e-04} & \texttt{        9} & \texttt{       25} & \texttt{       30} &                                          \\
                 &                         &                         & \texttt{                                1e+00} & \texttt{       30} & \texttt{       19} & \texttt{       10} &                                          \\
                 &                         &                         & \texttt{                                1e+04} & \texttt{       10} & \texttt{        7} & \texttt{        7} &                                          \\
  \midrule
  \texttt{    2} & \texttt{         14739} & \texttt{           729} & \texttt{                                1e-04} & \texttt{       13} & \texttt{       27} & \texttt{       32} &                                          \\
                 &                         &                         & \texttt{                                1e+00} & \texttt{       32} & \texttt{       19} & \texttt{       11} &                                          \\
                 &                         &                         & \texttt{                                1e+04} & \texttt{       11} & \texttt{        9} & \texttt{        9} &                                          \\
  \midrule
  \texttt{    3} & \texttt{        107811} & \texttt{          4913} & \texttt{                                1e-04} & \texttt{       13} & \texttt{       31} & \texttt{       33} &                                          \\
                 &                         &                         & \texttt{                                1e+00} & \texttt{       33} & \texttt{       19} & \texttt{       11} &                                          \\
                 &                         &                         & \texttt{                                1e+04} & \texttt{       11} & \texttt{        9} & \texttt{        9} &                                          \\
  \bottomrule
\end{tabular}

	\caption{Iteration numbers for the optimal control problem for the Stokes equation (\cref{subsec:experiments_distributed_optimal_control_Stokes}) required to reach a relative reduction by $10^{-6}$ of the initial residual norm \eqref{eq:total_residual_norm}. The parameter $\beta$ is fixed to $1$~\si{\obj\per\meter\tothe{5}\second\squared}. Notice that $\dim(Q_1) = \dim(V_1)$ and $\dim(Q_2) = \dim(V_2)$ holds.}
	\label{tab:convergence_OptimalControl3DStokes}
\end{table}

\section{Conclusion}
\label{sec:conclusion}

In this paper we introduced a new paradigm which may help design and find effective scalings of block-diagonal preconditioners for symmetric saddle-point problems.
To this end, we take note of the physical units of the primal and dual spaces involved in the problem.
The simple yet effective idea of the proposed dimensionally consistent preconditioning is to ensure that the block-diagonal preconditioner likewise respects these physical units.
As a consequence, the quantity monitored by \minres, \ie, the sum of squared preconditioner-induced residual norms, becomes physically meaningful.

In a number of examples covering fluid flow, elastic solid body deformation as well as optimal control problems, we showed that dimensional consistency can be achieved through scaling of the preconditioner blocks involving a natural combination of the parameters already present in the respective problem. 
In these examples it turns out that, simultaneously with dimensional consistency, we achieve robustness of the eigenvalues of the preconditioned operator with respect to \emph{all} problem parameters.
This observation suggests that parameter robust preconditioning and dimensional consistency are strongly related.
We thus conjecture that the concept of dimensional consistency can significantly facilitate the search for efficient and parameter robust preconditioners.

A number of open problems offer themselves for future investiagation.
One concerns the extension of the concept of dimensional consistency to non-self-adjoint problems which are solved, \eg, by \gmres.
Another question of interest is the investigation of problems with non-constant coefficients, \eg, inhomogeneous Lamé parameters in \cref{subsec:elasticity}.
While the evaluation of the relevant constants as in \eqref{eq:elasticity_relevant_constants} and \eqref{eq:elasticity_beta_unscaled} may become more involved, we conjecture that dimensional consistency remains a desirable property for the preconditioner also in the case of non-constant coefficients.
Details are postponed to future work.

\printbibliography

\end{document}